\documentclass[11pt]{amsart}
\usepackage{graphicx,amssymb,upref, amscd}

\let\<\langle
\let\>\rangle

\let\uml\"

\newcommand\Zplus{\mathbb{Z}^{+}}
\newcommand\Z{\mathbb{Z}}

\newcommand\Q{\mathbb{Q}}

\newcommand\R{{\mathbb{R}}}

\newcommand{\eps}{\varepsilon}

\def\p{{\mathfrak{p}}}

\newenvironment{Gal}{{\rm Gal\,}}{}

\newenvironment{Disc}{{\rm Disc}}{}

\newenvironment{Aut}{{\rm Aut}}{}

\newtheorem{theorem}{Theorem}[section]
\newtheorem{definition}[theorem]{Definition}
\newtheorem{lemma}[theorem]{Lemma}
\newtheorem{proposition}[theorem]{Proposition}
\newtheorem{proposition-definition}[theorem]{Proposition-Definition}
\newtheorem{corollary}[theorem]{Corollary}
\newtheorem{conjecture}[theorem]{Conjecture}

\newtheorem{question}[theorem]{Question}

\theoremstyle{definition}
\newtheorem{example}[theorem]{Example}

\theoremstyle{remark}
\newtheorem{remark}[theorem]{Remark}

\title[Eventually stable quadratic polynomials]{Eventually stable quadratic polynomials over $\mathbb{Q}$}

\author[DeMark]{David DeMark}  
\address{School of Mathematics \\ 
University of Minnesota \\
206 Church Street SE \\
Minneapolis, MN 55455, USA} 
\email{demar180@umn.edu}  

\author[Hindes]{Wade Hindes}  
\address{Department of Mathematics \\
Texas State University \\
601 University Drive \\
San Marcos, TX 78666, USA} 
\email{wmh33@txstate.edu}  

\author[Jones]{Rafe Jones}  
\address{Department of Mathematics and Statistics \\
Carleton College \\
1 North College St \\
Northfield, MN 55057, USA} 
\email{rfjones@carleton.edu} 

\author[Misplon]{Moses Misplon}
\address{Department of Mathematics and Statistics \\
Carleton College \\
1 North College St \\
Northfield, MN 55057, USA} 
\email{mzrmisplon@gmail.com}

\author[Stoll]{Michael Stoll}
\address{Department of Mathematics \\
Universit\"at Bayreuth\\
95440 Bayreuth, Germany} 
\email{Michael.Stoll@uni-bayreuth.de}

\author[Stoneman]{Michael Stoneman}
\address{Department of Mathematics and Statistics \\
Carleton College \\
1 North College St \\
Northfield, MN 55057, USA} 
\email{mstoneman@google.com}



\begin{document} 
 
\begin{abstract}
We study the number of irreducible factors (over $\Q$) of the $n$th iterate of a polynomial of the form $f_r(x) = x^2 + r$ for $r \in \Q$. When the number of such factors is bounded independent of $n$, we call $f_r(x)$ \textit{eventually stable} (over $\Q$). Previous work of Hamblen, Jones, and Madhu \cite{zdc} shows that $f_r$ is eventually stable unless $r$ has the form $1/c$ for some $c \in \Z \setminus \{0,-1\}$, in which case existing methods break down. We study this family, and prove that several conditions on $c$ of various flavors imply that all iterates of $f_{1/c}$ are irreducible. We give an algorithm that checks the  latter property for all $c$ up to a large bound $B$ in time polynomial in $\log B$.
We find all $c$-values for which the third iterate of $f_{1/c}$ has at least four irreducible factors, and all $c$-values such that $f_{1/c}$ is irreducible but its third iterate has at least three irreducible factors. This last result requires finding all rational points on a genus-2 hyperelliptic curve for which the method of Chabauty and Coleman does not apply; we use the more recent variant known as elliptic Chabauty. Finally, we apply all these results to completely determine the number of irreducible factors of any iterate of $f_{1/c}$, for all $c$ with absolute value at most $10^9$.
\end{abstract}

\maketitle

\section{Introduction}

Let $K$ be a field with algebraic closure $\overline{K}$, $f \in K[x]$, and $\alpha \in K$. For $n \geq 0$, let $f^n(x)$ be the $n$th iterate of $f$ (we take $f^0(x) = x$), and  $f^{-n}(\alpha)$ the set $\{\beta \in \overline{K} : f^n(\beta) = \alpha\}$. When $f^n(x) - \alpha$ is separable over $K$ for each $n \geq 1$, the set $T_f(\alpha) := \bigsqcup_{n \geq 0} f^{-n}(\alpha)$ \label{Tdef} acquires the structure of a rooted tree (with root $\alpha$) if we assign edges according to the action of $f$.

A large body of recent work has focused on algebraic properties of $T_f(\alpha)$, particularly the natural action of $\Gal(\overline{K}/K)$ on $T_f(\alpha)$ by tree automorphisms, which yields a homomorphism $\Gal(\overline{K}/K) \to \Aut(T_f(\alpha))$ called the arboreal Galois representation associated to $(f, \alpha)$.
A central question is whether the image of this homomorphism must have finite index in $\Aut(T_f(\alpha))$ (see \cite{survey} for an overview of work on this and related questions).
In the present article we study factorizations
of polynomials of the form $f^n(x) - \alpha$ in the case where $\alpha = 0$.
\begin{definition}
Let $K$ be a field and $f \in K[x]$, and $\alpha \in K$. We say $(f,\alpha)$ is \textbf{eventually stable over $K$} if there exists a constant $C(f,\alpha)$ such that
the number of irreducible factors over $K$ of $f^n(x) - \alpha$ is at most $C(f,\alpha)$ for all $n \geq 1$.

We say that $f$ is eventually stable over $K$ if $(f,0)$ is eventually stable.
\end{definition}

Apart from its own interest, eventual stability has proven to be a key link in at least two recent proofs of finite-index results for certain arboreal representations \cite{bridytucker, unicrit}. This is perhaps surprising given that eventual stability is a priori much weaker than finite index of the arboreal representation -- the former only implies that the number of Galois orbits on $f^{-n}(\alpha)$ is bounded as $n$ grows, which is an easy consequence of the latter.
There are other applications of eventual stability as well; for instance, if $f \in \Q[x]$ is eventually stable over $\Q$, then a finiteness result holds for $S$-integer points in the backwards orbit of $0$ under $f$ (see \cite[Section 3]{evstab} and \cite{sookdeo}).

The paper \cite{evstab} provides an overview of eventual stability and related ideas. That article defines a notion of eventual stability for rational functions, gives several characterizations of eventual stability, and states some general conjectures on the subject, all of which remain wide open. For example, a special case of \cite[Conjecture 1.2]{evstab} is the following: if $f \in \Q[x]$ is a polynomial of degree $d \geq 2$ such that $0$ is not periodic under $f$ (i.e. $f^n(0) \neq 0$ for all $n \geq 1$), then $f$ is eventually stable over $\Q$. Theorems 1.3 and 1.7 of \cite{evstab} also provide some of the few reasonably general results currently available on eventual stability. The proofs rely on generalizations of the Eisenstein criterion, and crucially assume good reduction of the rational function at the prime in question.

In this article, we address some of the conjectures in \cite{evstab} in cases where Eisenstein-type methods break down. One of our main results is the following:
\begin{theorem} \label{bigcomp}
Let $K = \Q$ and $f_r(x) = x^2 + r$ with $r = 1/c$ for $c \in \Z \setminus \{0, -1\}$. If $|c| \leq 10^9$, then $f_r$ is eventually stable over $\Q$ and $C(f_r,0) \leq 4$. More precisely, Conjecture \ref{mainconj} below holds for all $c$ with $|c| \leq 10^9$.
\end{theorem}

The family $x^2 + (1/c), \, c \in \Z \setminus \{0, -1\}$, is particularly recalcitrant. Eventual stability in this family (with $\alpha = 0$) is conjectured in \cite[Conjecture 1.4]{evstab}, and it is the only obstacle to establishing eventual stability (with $\alpha = 0$) in the family $x^2 + r, r \in \Q$. This is because \cite[Theorem 1.7]{evstab} handles the case when there is a prime $p$ with $v_p(r) > 0$. Moreover, \cite[Theorem 1.3]{evstab} uses $p = 2$ to establish eventual stability for $x^2 + 1/c$ when $c$ is odd, but when $c$ is even $x^2 + 1/c$ has bad reduction at $p = 2$, and Eisenstein-type methods break down completely.

We turn to methods inspired by \cite{quaddiv}, in particular various amplifications of \cite[Proposition 4.2]{quaddiv} and  \cite[Theorem 2.2]{itconst}, which state that the irreducibility of iterates can be proven by showing a certain sequence contains no squares. We prove the following theorem, which plays a substantial role in the proof of Theorem \ref{bigcomp}. For the rest of the article, we establish the following conventions:
$$
 \text{all irreducibility statements are over $\Q$} \qquad \text{and} \qquad \text{$r = 1/c$, where $c$ is a non-zero integer.}
$$
Also, we denote by $\Z \setminus \Z^2$ the set of integers that are not integer squares.

\begin{theorem} \label{main}
Let $f_r(x) = x^2 + r$ with $r = 1/c$ for $c \in \Z \setminus \{0, -1\}$. Then $f_r^n(x)$ is irreducible for all $n \geq 1$ if $c$ satisfies one of the following conditions:
\begin{enumerate}
\item $-c \in \Z \setminus \Z^2$ and $c < 0$;
\item $-c, c+1 \in \Z \setminus \Z^2$ and $c \equiv -1 \bmod{p}$ for a prime $p \equiv 3 \bmod{4}$;
\item $-c, c+1 \in \Z \setminus \Z^2$ and $c$ satisfies one of the congruences in Proposition \ref{congruences} (see Table 1).
\item $-c \in \Z \setminus \Z^2$ and $c$ is odd;
\item $-c \in \Z \setminus \Z^2$, $c$ is not of the form $4m^2(m^2-1), m \in \Z$, and

\begin{equation*}
\frac{\prod\limits_{p : 2 \nmid v_p(c)} p^{v_p(c)}}{\prod\limits_{p : 2 \mid v_p(c)} p^{v_p(c)}}
       > \frac{1.15}{|c|^{1/30}} .
\end{equation*}
This holds whenever $c$ is squarefree.

\item $c = k^2$ for some $k \geq 2$ and

\begin{equation*}
\frac{\prod\limits_{p : p \not\equiv 1 \bmod 4} p^{v_p(c)}}{\prod\limits_{p : p \equiv 1 \bmod 4} p^{v_p(c)}}
       > \frac{1.15}{|c|^{1/30}} .
\end{equation*}

\item  $c$ is not of the form $4m^2(m^2-1)$ with $m \in \Z$ and $1 \le c \le 10^{1000}$.
\end{enumerate}
\end{theorem}

Our next result  gives an explicit and relatively small bound for~$m$ such that
the irreducibility of~$f^m$ implies the irreducibility of all~$f^n$ (see Corollary~\ref{C:improve5.4}).
In the following, $\eps(c)$ is a function bounded above by $4$ and decreasing monotonically to $2$ as $c$ grows; for a precise definition, see p. \pageref{epsdef} in Section~\ref{mainresults}.

\begin{theorem} \label{stabtestthm}
  Let $f_r(x) = x^2 + r$ with $r = 1/c$ for $c \in \Z$ with $c \ge 4$.
  If $f^m$ is irreducible for
  \[ m = 1 + \Bigl\lfloor \log_2\Bigl(1 + \frac{\log 4 + \eps(c)/\sqrt{c}}%
                                                {\log(1 + 1/\sqrt{c})}\Bigr) \Bigr\rfloor ,
  \]
  then all $f^n$ are irreducible.
\end{theorem}

The methods used in the proof of this result can be used to derive a very efficient
algorithm that checks the condition in Theorem~\ref{stabtestthm} for all $c$ up to
a very large bound. This leads to a proof of case~(7) of Theorem~\ref{main},
which at the same time verifies Conjecture~\ref{stabconj}
below for all $c$ with $|c| \le 10^{1000}$. This is explained in Section~\ref{S:algo}.

We also prove results on unusual factorizations of small iterates in the family $x^2 + 1/c$.
 \begin{theorem} \label{maincases}
Let $f_r(x) = x^2 + r$ with $r = 1/c$ for $c \in \Z \setminus \{0, -1\}$, and let $k_n$ denote the number of irreducible factors of $f_r^n(x)$. Then
 \begin{enumerate}
 \item[(a)] We have $k_1 = k_2 = 2$ and $k_3 = 3$ if and only if $c = -16$. In this case $k_n = 3$ for all $n \geq 3$.
 \item[(b)] We have $k_3 \geq 4$ if and only if $c = -(s^2-1)^2$ for $s \in \{3, 5, 56\}$. In this case, $k_1 = 2$, $k_2 = 3$, and $k_n = 4$ for all $n \geq 3$.
 \item[(c)] We have $k_1 = 1$ and $k_3 \geq 3$ if and only if $c = 48$. In this case, $k_2 = 2$ and $k_n = 3$ for all $n \geq 3$.
 \end{enumerate}
 \end{theorem}

Observe that part (b) of Theorem \ref{maincases} shows that the bound $C(f_r,0) \leq 4$ in Theorem \ref{bigcomp} (and also Conjecture \ref{mainconj}) cannot be improved. Moreover, a consequence of Conjecture \ref{mainconj} is that $C(f_r,0) = 4$ if and only if $c = -(s^2-1)^2$ for $s \in \{3, 5, 56\}$.

To prove Theorem \ref{maincases}, we reduce the problem to finding all integer square values of certain polynomials (see Lemma \ref{squareirred} in Section \ref{redsec} for details.) The curve that arises in this way in the proof of part (c) of Theorem \ref{maincases} is of particular interest, as it is a hyperelliptic curve of genus two, whose Jacobian has rank two:
\begin{equation} \label{maincurve}
y^2 = 8x^6 - 12x^4 - 4x^3 + 4x^2 + 4x + 1.
\end{equation}
Because the genus and the rank of the Jacobian coincide, we cannot apply the well-known method of Chabauty and Coleman. On the other hand, we are able to use a variant of the standard method, called elliptic Chabauty \cite{Bruin, Flynn}, to prove:
\begin{theorem} \label{ellchab}
The only rational points on the curve \eqref{maincurve} are those with $x \in \{-2, -1, 0, 1\}$.
\end{theorem}
The idea is the following. Given an elliptic curve $E$ over a number field $K$ and a map $\phi:E\rightarrow \mathbb{P}^1$, then one can often compute the set of points in $E(K)$ mapping to $\mathbb{P}^1(\mathbb{Q})$ as long as the rank of $E(K)$ is strictly less than the degree of the extension $K/\mathbb{Q}$; this method is known as elliptic Chabauty. Moreover, in certain situations, one can use a combination of descent techniques and elliptic Chabauty to determine the full set of rational points on a curve $C$ (of higher genus) defined over $\mathbb{Q}$; see, for instance, the proof of Theorem \ref{ellchab}. Moreover, under suitable conditions, several components of the elliptic Chabauty method are implemented in MAGMA~\cite{magma}, and we make use of these implementations here. Our code verifying the calculations in the proof of Theorem 1.5 can be found within the file called \emph{Elliptic Chabauty} at:
${\tt{https://sites.google.com/a/alumni.brown.edu/whindes/research}}\vspace{.05cm}$.

The above results furnish evidence for several conjectures. The first is a refinement of Conjecture 1.4 of \cite{evstab}, which states that $x^2 + \frac{1}{c}$ is eventually stable for $c \in \Z \setminus \{0,-1\}$.
 \begin{conjecture} \label{mainconj}
Let $f_r(x) = x^2 + r$ with $r = 1/c$ for $c \in \Z \setminus \{0, -1\}$. Then $f_r$ is eventually stable and $C(f_r, 0) \leq 4$.
More precisely, let $k_n$ denote the number of irreducible factors of $f_r^n(x)$. Then
\begin{enumerate}
\item If $c = -m^2$ for $m > 0$ with $m + 1 \in \Z \setminus \Z^2$ and $m \neq 4$, then $k_n = 2$ for all $n \geq 1$.
\item If $c = -16$, then $k_1 = k_2 = 2$ and $k_n = 3$ for all $n \geq 3$.
\item If $c = -(s^2-1)^2$ for $s \in \Z \setminus \{3, 5, 56\}$, then $k_1 = 2$ and $k_n = 3$ for all $n \geq 2$.
\item If $c = -(s^2-1)^2$ for $s \in \{3, 5, 56\}$, then $k_1 = 2$, $k_2 = 3$, and $k_n = 4$ for all $n \geq 3$.
\item If $c = 4m^2(m^2-1)$ for $m \in \Z$, $m \geq 3$, then $k_1 = 1$ and $k_n = 2$ for all $n \geq 2$.
\item If $c = 48$, then $k_1 = 1$, $k_2 = 2$, and $k_n = 3$ for all $n \geq 3$.
\item If $c$ is not in any of the above cases, then $k_n = 1$ for all $n \geq 1$.
\end{enumerate}
\end{conjecture}

We remark that case (7) of Conjecture \ref{mainconj} is precisely the case where $f_r^2(x)$ is irreducible (see Proposition \ref{irredprop}) and thus case (7) asserts that if $f_r^2(x)$ is irreducible, then $f_r^n(x)$ is irreducible for all $n \geq 1$. We state this as its own conjecture:
\begin{conjecture} \label{stabconj}
Let $f_r(x) = x^2 + r$ with $r = 1/c$ for $c \in \Z \setminus \{0, -1\}$. If $f_r^2(x)$ is irreducible, then $f_r^n(x)$ is irreducible for all $n \geq 1$.
\end{conjecture}

 As mentioned above, we have verified this conjecture for all $c$ with $|c| \le 10^{1000}$.

Observe that Conjecture \ref{mainconj} gives a uniform bound for $k_n$, in contrast to Conjecture 1.4 of \cite{evstab}. It would be of great interest to have a similar uniform bound for $f_r(x)$ as $r$ is allowed to vary over the entire set $\Q \setminus \{0, -1\}$ (as opposed to just the reciprocals of integers, as in Conjecture \ref{mainconj}). We pose here a much more general question. Given a field $K$, call $f \in K[x]$ \textit{normalized} (the terminology \textit{depressed} is also sometimes used, especially for cubics) if $\deg f = d \geq 2$ and $f(x) = a_dx^d + a_{d-2}x^{d-2} + a_{d-3}x^{d-3} + a_1x + a_0$. Note that every $f \in K[x]$ of degree not divisible by the characteristic of $K$ is linearly conjugate over $K$ to a normalized polynomial.

\begin{question} \label{unifbound}
Let $K$ be a number field and fix $d \geq 2$. Is there a constant $\kappa$ depending only on $d$ and $[K : \Q]$ such that, for all normalized $f \in K[x]$ of degree $d$ such that $0$ is not periodic under $f$, and all $n \geq 1$, $f^n(x)$ has at most $\kappa$ irreducible factors? In the case where $K = \Q$, $d = 2$, and $f$ is taken to be monic, does the same conclusion hold with $\kappa = 4$?
\end{question}

It is interesting to compare Question \ref{unifbound} to \cite[Question 19.5]{trends}, where a similar uniform bound is requested, but under the condition that $f^{-1}(0) \cap \mathbb{P}^1(K) = \emptyset$.

We close this introduction with some further comments on our methods, as well as the statement of one additional result (Theorem \ref{densitythm}) on the densities of primes dividing orbits of polynomials of the form $x^2 + 1/c$.

A primary tool in our arguments is the following special case of \cite[Theorem 2.2]{itconst}: for $n \geq 2$, $f_r^n$ is irreducible provided that $f_r^{n-1}$ is irreducible and $f_r^n(0)$ is not a square in $\Q$. The proof of this relies heavily on the fact that $f_r$ has degree 2, and is essentially an application of the multiplicativity of the norm map.
 Using ideas from \cite[Theorem 2.3 and discussion preceding]{itconst}, one obtains the useful amplification (proven in Section \ref{preliminaries}) that for $n \geq 2$, $f_r^n$ is irreducible provided that $f_r^{n-1}$ is irreducible and neither of $(f_r^{n-1}(0) \pm \sqrt{f_r^n(0)})/2$ is a square in $\Q$.
When $r = 1/c$, we have
$$f_r(0) = 1/c, \quad f_r^2(0) = (c + 1)/c^2, \quad f_r^3(0) = (c^3 + c^2 + 2c + 1)/c^4,$$
and so on. The numerator of $f_r^n(0)$ is obtained by squaring the numerator of $f_r^{n-1}(0)$, and adding $c^{2^{n-1}-1}$.
We thus introduce the family of sequences
\begin{equation} \label{andef}
a_1(c) = 1, \qquad a_{n}(c) = a_{n-1}(c)^2 + c^{2^{n-1} - 1} \quad \text{for $n \geq 2$}.
\end{equation}
To ease notation, we often suppress the dependence on $c$, and write $a_1, a_2,$ etc. We can then translate the results of the previous paragraph to:
\begin{lemma} \label{biglemma}
Suppose that $c \in \Z \setminus \{0\}$, $r = 1/c$, and $f_r^2$ is irreducible. Let $a_n = a_n(c)$ be defined as in \eqref{andef}, and set
\begin{equation} \label{bndef}
b_n := \frac{a_{n-1} + \sqrt{a_{n}}}{2} \in \overline{\mathbb{Q}}.
\end{equation}
 If for every $n \geq 3$, $b_n$ is not a square in $\Q$ (which holds in particular if $a_n$ is not a square in $\Q$), then $f_r^n(x)$ is irreducible for all $n \geq 1$.
\end{lemma}

We make the following conjecture, which by Lemma \ref{biglemma} immediately implies Conjecture \ref{stabconj}:
\begin{conjecture} \label{nosquare}
Let $b_n = b_n(c)$ be defined as in \eqref{bndef}. If $c \in \Z \setminus \{0, -1\}$, then $b_n$ is not a square in $\mathbb{Q}$ for all $n \geq 3$.
\end{conjecture}

Conjecture \ref{nosquare} also has strong implications for the density of primes dividing orbits of $f_r$.
We define the orbit of $t \in \Q$ under $f_r$ to be the set $O_{f_r}(t) = \{t, f_r(t), f_r^2(t), \ldots\}$, and we say that a prime $p$ divides $O_{f_r}(t)$ if there is at least one non-zero $y \in O_{f_r}(t)$ with $v_p(y) > 0$. The natural density of a set $S$ of prime numbers is defined to be
$$D(S) = \lim_{B \to \infty} \frac{\#\{p \leq B : p \in S\}}{\#\{p \leq B\}}.$$
Note that the elements of $O_{f_r}(t)$ also form a nonlinear recurrence sequence, where the relation is given by application of $f_r$. The problem of finding the density of prime divisors in recurrences has an extensive literature in the case of a linear recurrence; see the discussion and brief literature review in \cite[Introduction]{quaddiv}. The case of non-linear recurrences is much less-studied, though there are some recent results \cite{zdc, quaddiv, looper}. The following theorem is an application of
\cite[Theorem 1.1, part~(2)]{zdc}.

\begin{theorem} \label{densitythm}
Let $c \in \Z$, let $r = 1/c$, suppose that $-c$ and $c+1$ are non-squares in $\Q$, and assume that Conjecture \ref{nosquare} holds for $c$. Then
\begin{equation} \label{denseq0}
\text{for any $t \in \Q$ we have $D(\{\text{$p$ prime : $p$ divides $O_{f_r}(t)$\}}) = 0$. }
\end{equation}
\end{theorem}
We remark that in each of the cases of Theorem \ref{main}, we show that Conjecture \ref{nosquare} holds for $c$. Hence in cases (2), (3), and (6) of Theorem \ref{main}, and also in cases (1), (4), and (5) with the additional hypothesis that $c+1$ is not a square in $\Q$, we have that \eqref{denseq0} holds. We also note that when the hypotheses of Theorem \ref{densitythm} are satisfied, we obtain certain information on the action of $G_\Q$ on $T_f(0) = \bigsqcup_{n \geq 0} f^{-n}(0)$; see Section \ref{density}.

A complete proof of Conjecture \ref{nosquare} appears out of reach at present. One natural approach is to prove the stronger statement that $a_n$ is not a square for each $n \geq 3$, or equivalently that
the curve
\begin{equation} \label{curves}
C_n : y^2 = a_n(c)
\end{equation}
has no integral points with $c \not\in \{0, -1\}$ for any $n \geq 3$. It is easy to see that $a_n(c)$ is separable as a polynomial in $c$ (one considers it as a polynomial in $\Z/2\Z[c]$, where it is relatively prime to its derivative), and because the degree of $a_n(c)$ is $2^{n-1}-1$, it follows from standard facts about hyperelliptic curves that the genus of $C_n$ is $2^{n-2} - 1$. Siegel's theorem then implies that there are only finitely many $c$ with $a_n(c)$ a square for given $n \geq 3$. However, the size of the genus of $C_n$ prevents us from explicitly excluding the presence of integer points save in the cases of $n = 3$ and $n = 4$ (see Proposition \ref{a3nonsquare}).
One idea that has been used to show families of integer non-linear recurrences contain no squares (see e.g. \cite[Corollary 1.3]{stoll} or \cite[Lemma 4.3]{quaddiv}) is to show that sufficiently large terms of each sequence are sandwiched between squares: they are generated by adding a small number to a large square. In the case of the family $a_n(c)$, however, the addition of the very large term $c^{2^{n-1}-1}$ to the square $a_{n-1}^2$ ruins this approach. A similar problem is encountered in a family of important two-variable non-linear recurrence sequences first considered in \cite{galrat} (see \cite[Theorem 1.8]{galrat}). The main idea used in \cite{galrat} to show the recurrence contains no squares is to rule out certain cases via congruence arguments. This is the essence of our method of proof for cases (2) and (3) of Theorem \ref{main}. Subsequently Swaminathan \cite[Section 4]{swaminathan} amplified these congruence arguments and gave new partial results using the idea of sandwiching terms of the sequence between squares. In the end each of these methods succeeds in giving only partial results, applicable to $c$-values satisfying certain arithmetic criteria. It would be of great interest to have a proof of Conjecture \ref{nosquare} for $c$-values satisfying some analytic criterion, e.g., for all $c$ sufficiently large. Case (1) of Theorem \ref{main} provides one result with this flavor, but at present no other results are known.

\textbf{Acknowledgements:} We thank Jennifer Balakrishnan for conversations related to the proof of Theorem \ref{ellchab}, and the anonymous referee for many helpful suggestions.

\section{The case where $f_r(x)$ or $f_r^2(x)$ is reducible} \label{redsec}

We begin by studying the factorizations of iterates of $f_r(x)$ when either $f_r(x)$ or $f^2_r(x)$ is reducible. The behavior of higher iterates becomes harder to control because of the presence of multiple irreducible factors of the first two iterates, but we are still able to give some results. At the end of this section we prove Theorem \ref{maincases}, which gives a complete characterization of certain subcases.

\begin{proposition} \label{irredprop}
Let $f_r(x) = x^2 + r$ with $r = 1/c$ for $c \in \Z \setminus \{0, -1\}$. Then $f_r(x)$ is reducible if and only if $c = -m^2$ for $m \in \Z$. If $f_r(x)$ is irreducible, then $f_r^2(x)$ is reducible if and only if $c = 4m^2(m^2 - 1)$ for $m \in \Z$.
\end{proposition}

\begin{proof}
The first statement is clear. Assume now that $f_r(x)$ is irreducible over $\Q$. Let $\alpha$ be a root of $f_r^2(x)$, and observe that $f_r(\alpha)$ is a root of $f_r(x)$, and by the irreducibility of $f_r(x)$, we have $[\Q(f_r(\alpha)) : \Q] = 2$. Now $f_r^2(x)$ is irreducible if and only if $[\Q(\alpha) : \Q] = 4$, which is equivalent to $[\Q(\alpha) : \Q(f_r(\alpha))] = 2$.
But $\alpha$ is a root of $f_r(x) - f_r(\alpha) = x^2 + r -f_r(\alpha)$, and so $[\Q(\alpha) : \Q(f_r(\alpha))] = 2$ is equivalent to $f_r(\alpha) - r$ not being a square in $\Q(f_r(\alpha))$.

Without loss of generality, say $f_r(\alpha) = \sqrt{-r}$. Then $f_r(\alpha) - r$ is a square in $\Q(f_r(\alpha))$ if and only if there are $s_1, s_2 \in \Q$ with
$$-r + \sqrt{-r} = (s_1 + s_2\sqrt{-r})^2 =  s_1^2 - rs_2^2 + 2s_1s_2\sqrt{-r}.$$
This holds if and only if $2s_1s_2 = 1$ and $s_1^2 - rs_2^2 = -r$. Substituting $s_2 = 1/(2s_1)$ into the second equation and multiplying through by $s_1^2$ gives $s_1^4 + rs_1^2 - r/4 = 0$, which by the quadratic formula holds if and only if
\begin{equation} \label{square}
s_1^2 = \frac{-r \pm \sqrt{r^2 + r}}{2}
\end{equation}
or equivalently, $2c(-1 \pm \sqrt{1+c})$ is an integer square (here we have written $1/c$ for $r$ and multiplied both sides of \eqref{square} by $4c^2$). If $c < -1$, then $\sqrt{1 + c}$ is irrational, so we may assume $c > 0$. We may then discard the $-$ part of the $\pm$, since integer squares are positive. Writing $c =k^2-1$ for $k > 0$, we then obtain that $2(k^2-1)(-1 + k) = 2(k+1)(k-1)^2$ is a square, whence $k+1 = 2m^2$ for some integer $m$. Thus $c = k^2 - 1 = (2m^2-1)^2 - 1 = 4m^4 - 4m^2$, as desired.
\end{proof}

We now give a lemma, closely related to \cite[Proposition 4.2]{quaddiv}, which we will use often in the sequel.

\begin{lemma} \label{squareirred}
Let $K$ be a field of characteristic not equal to 2, let $g \in K[x]$ be monic of degree $d \geq 1$ and irreducible over $K$, and let $f(x)$ be monic and quadratic with critical point $\gamma$. If no element of $$\{(-1)^dg(f(\gamma))\} \cup \{g(f^n(\gamma)) : n \geq 2\}$$ is a square in $K$, then $g(f^n(x))$ is irreducible over $K$ for all $n \geq 1$.
\end{lemma}

\begin{proof}
Let $f(x) = x^2 + bx + c$, so that $\gamma = -b/2$. We proceed by induction on $n$, with the $n = 0$ case covered by the irreducibility of $g(x)$. Assume then that $g(f^{n-1}(x))$ is irreducible over $K$ for some $n \geq 1$, and let $d_1$ be the degree of $g(f^{n-1}(x))$. By Capelli's Lemma (\cite[Lemma 0.1]{fein}), $g(f^n(x))$ is irreducible over $K$ if and only if for any root $\beta$ of $g(f^{n-1}(x))$,
we have $f(x) - \beta$ is irreducible over $K(\beta)$, or equivalently (because $K$ has characteristic different from 2),
$\Disc (f(x) - \beta) = b^2 - 4c + 4\beta$ is not a square in $K(\beta)$.

This must hold if $N_{K(\beta)/K} (b^2 - 4c + 4\beta)$ is not a square in $K$.  The Galois conjugates of $b^2 - 4c + 4\beta$ are precisely $b^2 - 4c + 4\alpha$ as $\alpha$ varies over all roots of $g(f^{n-1}(x))$. Thus
\begin{eqnarray*}
N_{K(\beta)/K} (b^2 - 4c + 4\beta) & = & (-4)^{d_1} \prod_{\text{$\alpha$ root of $g \circ f^{n-1}$}}
\left[ \left(- \frac{b^2}{4} + c \right) - \alpha \right]\\
& = & (-4)^{d_1} \cdot g(f^{n-1}(-b^2/4 + c)) \\ 
& = & (-4)^{d_1} \cdot g(f^{n-1}(f(\gamma))),
\end{eqnarray*}
where the second equality holds because $g(f^{n-1}(x))$ is monic. Now $d_1$ is odd if and only if $n = 1$ and $d$ is odd, which proves the Lemma.
\end{proof}

\subsection{The case of $f_r$ reducible}
When $c = -m^2$ for some $m \geq 1$, we fix the notation
\begin{equation} \label{gnotation}
g_1(x) = x - \frac{1}{m}\qquad \text{and} \qquad g_2(x) = x + \frac{1}{m},
\end{equation}
so that $f_r(x) = g_1(x)g_2(x)$. We exclude the case $m = 1$ in what follows, as in that case $f_r(x)$ is not eventually stable (see \cite[discussion following Corollary 1.5]{evstab}).

\begin{proposition} \label{firstcaseprop}
Let $r = 1/c$ and $c = -m^2$ for $m \geq 2$. Let $g_1$ and $g_2$ be as in \eqref{gnotation}. Then the following hold.
\begin{enumerate}
\item We have $g_2(f_r(x))$ irreducible, while $g_1(f_r(x))$ factors if and only if $m+1$ is a square in $\Q$.
\item If $g_1(f_r(x))$ is irreducible, then $g_1(f_r^n(x))$ is irreducible for all $n \geq 2$.
\item If every term of the sequence $\{g_2(f_r^{i}(0))\}_{i\geq 2}$ is a non-square in $\Q$, then
$g_2(f_r^n(x))$ is irreducible for all $n \geq 2$.
\end{enumerate}
\end{proposition}

\begin{proof}
The first item follows from observing that
$g_1(f_r(x)) = x^2-\frac{m+1}{m^2}$ and $g_2(f_r(x)) = x^2+\frac{m-1}{m^2}$. The latter is irreducible because $m \geq 2$ implies $(m-1)/m^2 > 0$.
Item (3) is an immediate consequence of item (1) and Lemma \ref{squareirred} (with $g = g_2 \circ f_r$ and $f = f_r$). To prove item (2), observe that
$g_1(f_r^n(0)) = f_r^n(0) - \frac{1}{m}.$
However, one easily checks that $x^2 - \frac{1}{m^2}$ maps the interval $(-1/m, 0)$ into itself, and in particular, $f_r^n(0) < 0$ for all $n \geq 1$. Thus $g_1(f_r^n(0)) < 0$ as well, and hence cannot be a square in $\Q$. Item (2) now follows from Lemma \ref{squareirred} with $g = g_1 \circ f_r$ and $f = f_r$.
\end{proof}

\begin{proposition} \label{m4prop}
Let $r = 1/c$ and $c = -m^2$ for $m \geq 2$, and let $g_1$ and $g_2$ be as in \eqref{gnotation}. Then $g_2(f_r^2(0))$ is a square in $\Q$ if and only if $m = 4$. Moreover, $g_2(f_r^2(x))$ is reducible if and only if $m = 4$.
\end{proposition}

\begin{proof}
Observe that
\begin{equation*} 
g_2(f_r^2(0)) = \frac{m^3 - m^2 + 1}{m^4},
\end{equation*}
and hence $g_2(f_r^2(0))$ is a square in $\Q$ if and only if the elliptic curve $y^2 = x^3 - x^2 + 1$ has an integral point with $x = m$. This is curve 184.a1 in the LMFDB \cite{lmfdb}, and has only the integral points
$(0,\pm 1), (1,\pm 1), (4,\pm 7)$. Because $m \geq 2$, the only $m$-value for which $g_2(f_r^2(0))$ is a square is $m = 4$.
If $m \neq 4$, then \cite[Proposition 4.2]{quaddiv} (or the proof of Lemma \ref{squareirred}, with $g = g_2 \circ f_r$ and $f = f_r$) shows that $g_2(f_r^2(x))$ is irreducible. On the other hand, if $m = 4$, then
\begin{equation*}
\label{m4fact}
g_2(f_r^2(x)) = (x^2 - x + 7/16)(x^2 + x + 7/16), 
\end{equation*}
showing that $g_2(f_r^2(x))$ is reducible. We return to the analysis of the case $m = 4$ in Proposition \ref{m4}.
\end{proof}

We now seek to give congruence conditions on $m$ that ensure the sequence $(g_2(f_r^n(0)))_{n \geq 2}$ contains no squares in $\Q$. Prime factors of the numerators of the terms of this sequence are often related to each other. To formalize this, we require the following definition.

\begin{definition} \label{rigid}
A sequence $(s_n)_{n \geq 1}$ is a \textit{rigid divisibility sequence} if for all primes $p$ we have the following:
\begin{enumerate}
\item if $v_p(s_n) = e > 0$, then $v_p(s_{mn}) = e$ for all $m \geq 1$, and
\item if $v_p(s_n) > 0$ and $v_p(s_j) > 0$, then $v_p(s_{\gcd(n,j)}) > 0$.
\end{enumerate}
\end{definition}
\begin{remark}
If $(s_n)_{n \geq 1}$ is a rigid divisibility sequence and $s_1 = 1$, then from (2) it follows that if $p \mid \gcd(s_n, s_{n-1})$ then $p \mid s_1 = 1$, which is impossible. Hence
$\gcd(s_n, s_{n-1}) = 1$ for all $n \geq 2$. A similar argument shows that for $q$ prime we have $\gcd(s_q, s_i) = 1$ for all $1 \leq i < q$.
\end{remark}

\begin{proposition} \label{congruences2}
Let $r = 1/c$ and $c = -m^2$ for $m \geq 2$, and let $g_2$ be as in \eqref{gnotation}.
Then $g_2(f_r^n(x))$ is irreducible for all $n \geq 2$ provided that $m \neq 4$ and at least one of the following holds:
\begin{align*}
m  & \equiv  3 &  \pmod{4} & \qquad & m  &\equiv 3 & \pmod{5} \\
m &  \equiv  2, 5, 6 & \pmod{7} & \qquad & m & \equiv 4, 6, 7  & \pmod{11} \\
m &  \equiv  8, 10 & \pmod{13} & \qquad & m & \equiv 2, 4, 7, 8, 9, 11, 15  & \pmod{17}  \\
m &  \equiv  3, 5, 11 & \pmod{19} & \qquad & m & \equiv 9, 11, 14, 15, 18, 20, 21, 22  & \pmod{23} \\
m &  \equiv  3, 19, 26 & \pmod{29} & \qquad & m & \equiv 2, 12, 30  & \pmod{31} \\
m &  \equiv  6, 20 & \pmod{37} & \qquad & m & \equiv 12, 14, 27, 29  & \pmod{41} \\
m &  \equiv  15, 21, 30 & \pmod{43} & \qquad & m & \equiv 9, 22, 38, 46  & \pmod{47}
\end{align*}
If in addition $m-1$ is not a square in $\Q$, then the following congruences also suffice:
\begin{align*}
m & \equiv 2 & \pmod{3} & \qquad & m  &\equiv 5 & \pmod{8} \\
m  &\equiv 10 & \pmod{11} & \qquad & m & \equiv 18 & \pmod{19} \\
 m  &\equiv 2, 13 & \pmod{23} & \qquad & m & \equiv 8, 10, 14 & \pmod{29}\\
 m  &\equiv 9, 26 & \pmod{31} & \qquad & m & \equiv 13, 31 & \pmod{37}\\
m  &\equiv 3, 11, 19, 37, 38 & \pmod{41} & \qquad & m & \equiv 22, 36, 39, 42 & \pmod{43} \\
m  &\equiv 3, 10 & \pmod{47} \\
\end{align*}
\end{proposition}

\begin{proof}
By part (3) of Proposition \ref{firstcaseprop}, it suffices to show that $g_2(f_r^n(0))$ is not a square in $\Q$ for all $n \geq 2$. Note that for each $n \geq 1$, $g_2(f_r^{n-1}(0))$ is a positive rational number with denominator $m^{2^{n-1}}$, and numerator prime to $m$. We take $w_n$ to be the numerator of $g_2(f_r^{n-1}(0))$. We first observe that the proof of \cite[Proposition 5.4]{quaddiv} shows that the sequence $(w_n)_{n \geq 1}$ is a rigid divisibility sequence. In particular, if $w_2$ is not a square in $\Q$, then because $w_2 > 0$ we must have some prime $p$ dividing $w_2$ to odd multiplicity, and the rigid divisibility condition implies that $w_{2j}$ is not a square for all $j \geq 2$. A similar argument shows that if $w_3$ is not a square in $\Q$, then neither is $w_{3j}$ for all $j \geq 1$.

By Proposition \ref{m4prop} and our assumption that $m \neq 4$, we have that $g_2(f_r^2(0))$ is not a square in $\Q$. It follows that $w_{3j}$ is a non-square for all $j \geq 1$.

Now for a given modulus $k$ and $m \not\equiv 0 \bmod{k}$, the sequence $(g_2(f_r^n(0)) \bmod{k})_{n \geq 1}$ eventually lands in a repeating cycle, and we search for values of $k$ and congruences classes of $m$ modulo $k$ such that $g_2(f_r^n(0)) \bmod{k}$ fails to be a square for all $n \geq 2$. Note that this method works even when $g_2(f_r^{3j-1}(0)) \bmod{k}$ is a square for all $j \geq 1$, since we have shown in the previous paragraph that $w_{3j}$ is a non-square for all $j \geq 1$. A computer search yields the congruences given in the first part of the proposition. If in addition $m-1$ is a non-square in $\Q$, then we have $w_{2j}$ not a square in $\Q$ for all $j \geq 1$, and the congruences in the second part of the proposition show that $w_{2j+1} = g_2(f_r^{2j}(0)) \bmod{k}$ is a non-square for all $j \geq 1$.
\end{proof}

\begin{proposition} \label{congruences21}
Let $r = 1/c$ and $c = -m^2$ for $m \geq 2$, and let $g_2$ be as in \eqref{gnotation}. If $m \equiv -1 \bmod{p}$ for a prime $p \equiv 7 \bmod{8}$, then $g_2(f_r^n(x))$ is irreducible for all $n \geq 2$.  The same conclusion holds if  $m-1$ is not a square in $\Q$ and $m \equiv -1 \bmod{p}$ for a prime $p \equiv 3 \bmod{8}$.
\end{proposition}

\begin{proof}
By part (3) of Proposition \ref{firstcaseprop}, it suffices to show that $g_2(f_r^n(0))$ is not a square in $\Q$ for all $n \geq 2$. We have $c = -m^2 \equiv -1 \bmod{p}$, and so $(f_r^n(0) \bmod{p})_{n \geq 0}$ is the sequence $0, -1, 0, -1 \ldots$. Thus
$(g_2(f_r^n(0)) \bmod{p})_{n \geq 0}$ is the sequence $-1, -2, -1, -2, -1, \ldots$. If $p \equiv 7 \bmod{8}$, then both $-1$ and $-2$ are non-squares modulo $p$, and the proof is complete. If $p \equiv 3 \bmod{8}$, then $-1$ is a non-square modulo $p$ but $-2$ is a square, meaning we can only conclude that $g_2(f_r^{2j}(0))$ is a non-square in $\Q$ for $j \geq 1$. However, as in the proof of Proposition \ref{congruences2}, this implies that $b_{2j+1}$ is a non-square for all $j \geq 1$. If in addition $m-1$ is not a square, then $b_{2j}$ is not a square for all $j \geq 1$, completing the proof.
\end{proof}

Propositions \ref{congruences2} and \ref{congruences21} allow us to prove a case of Theorem \ref{bigcomp}. Recall that $g_1(f_r^n(x))g_2(f_r^n(x)) = f_r^{n+1}(x)$.

\begin{corollary} \label{g2cor}
Let $r = 1/c$ and $c = -m^2$ for $m \geq 2$, and let $g_2$ be as in \eqref{gnotation}. Suppose that $m \neq 4$ and $m^2 \leq 10^9$. Then $g_2(f_r^n(x))$ is irreducible for all $n \geq 1$. If in addition $m+1$ is not a square in $\Q$, then $f_r^n(x)$ is a product of two irreducible factors for all $n \geq 1$.
\end{corollary}

\begin{proof}
By part (3) of Proposition \ref{firstcaseprop}, it suffices to show that $g_2(f_r^n(0))$ is not a square in $\Q$ for all $n \geq 2$. Because $m \neq 4$, we may apply both Propositions \ref{congruences2} and \ref{congruences21}.
The first group of congruences in Proposition \ref{congruences2} applies to all $m$ with $2 \leq m \leq 10^{9/2}$ except for a set of 1326 $m$-values. After applying the first part of Proposition \ref{congruences21}, that number decreases to 1021. Of these, 13 have the property that $m-1$ is a square. We apply the second group of congruences in Proposition \ref{congruences2} and the second part of Proposition \ref{congruences21} to the remaining 1008 values, and only 196 survive. This leaves $209$ values of $m$ that we must handle via other methods.

To do this, we employ a new method to search for primes $p$ such that $g_2(f_r^n(0))$ is a non-square modulo $p$ for all but finitely many $n$. We search for $p$ such that:
\begin{multline} \label{goal}
\text{the sequence $(g_2(f_r^n(0)) \bmod{p})_{n \geq 0}$ eventually assumes} \\ \text{ a non-square constant value or eventually cycles between two distinct} \\ \text{values, both of which are non-squares modulo $p$.}
\end{multline}
If we find such a $p$, it implies that all but finitely many terms of the sequence $(g_2(f_r^n(0)))_{n \geq 2}$ are non-squares in $\Q$. We then reduce modulo other primes to show that the remaining terms are non-squares.

The method proves quite effective. Of the 209 $m$-values left over from the first paragraph of this proof, all have a prime $p < 500$ that satisfies \eqref{goal}. 
For each such $m$ and $p$, we take the finitely many terms of the sequence $(g_2(f_r^n(0)))_{n \geq 2}$ that have still not been proven non-square by \eqref{goal}, and reduce modulo small primes until all have been proven non-square.
The $m$-value producing the largest number of such terms is $m = 4284$, where we must check that each of $g_2(f_r(0)), g_2(f_r^2(0)), \ldots, g_2(f_r^{34}(0))$ is a non-square. In all cases the desired result is achieved by reducing modulo primes less than 100.
\end{proof}

We now consider the case $m = 4$. As shown in Proposition \ref{m4prop}, it is the only one with $m \geq 2$ for which $g_2(f_r^2(x))$ is reducible; indeed, we have
\begin{equation} \label{g2notation}
g_2(f_r^2(x)) = (x^2 - x + 7/16)(x^2 + x + 7/16) =: g_{21}(x)g_{22}(x),
\end{equation}
and we note that both $g_{21}(x)$ and $g_{22}(x)$ are irreducible.

\begin{proposition} \label{m4}
Let $r = -1/16$ and let $g_{21}$ and $g_{22}$ be as in \eqref{g2notation}. For all $n \geq 1$, both $g_{21}(f_r^n(x))$ and $g_{22}(f_r^n(x))$ are irreducible for all $n \geq 1$. In particular, $f_r^n(x)$ has precisely three irreducible factors for all $n \geq 3$.
\end{proposition}

\begin{proof}
Because $m+1$ is not a square, Proposition \ref{firstcaseprop} shows that $g_1(f_r^n(x))$ is irreducible for all $n \geq 1$. By Lemma \ref{squareirred} and the fact that $g_{21}$ and $g_{22}$ have even degree, it suffices to prove that neither $g_{21}(f_r^n(0))$ nor $g_{22}(f_r^n(0))$ is a square in $\Q$ for all $n \geq 1$.
Observe that $f_r^n(0) \equiv 5 \bmod{11}$ for $n \geq 3$, and $g_{21}(5) \equiv 6 \bmod{11}$. Because $6$ is a non-square modulo $11$, we must only verify that neither of $g_{21}(f_r(0))$ or
$g_{21}(f_r^2(0))$ is a square in $\Q$. The former is $129/256$ and the latter is $(19\cdot 1723)/2^{16}$, neither of which is a square in $\Q$. For $g_{22}(f_r^n(0))$ we have a simpler argument using $p = 5$: observe that $g_{22}(0) \equiv g_{22}(-1) \equiv 2 \bmod{5}$ and $f_r^n(0) \equiv 0$ or $-1 \bmod{5}$ for all $n \geq 1$.
\end{proof}

We now consider the case where $m+1$ is a square. Say $m + 1 = s^2$ with $s \geq 2$, so that $f_r(x) = x^2 - 1/m^2 = x^2 - 1/(s^2-1)^2$. We have
\begin{equation} \label{hdef}
g_1(f_r(x)) = x^2 - \frac{m+1}{m^2} = \left(x - \frac{s}{s^2 - 1} \right)\left(x + \frac{s}{s^2 - 1} \right) =: h_1(x)h_2(x).
\end{equation}
Now $h_1(f_r(x)) = x^2 - \frac{s^3 - s + 1}{(s^2 - 1)^2}$. Thus $h_1(f_r(x))$ is irreducible unless $s$ is the $x$-coordinate of an integral point on the elliptic curve $y^2 = x^3 - x + 1$. This is curve 92.a1 in LMFDB, and has an unusually large number of integral points: $(0, \pm 1), (1, \pm 1), (-1, \pm 1)(3, \pm 5), (5, \pm 11), (56, \pm 419)$. We assume for a moment that $s \not\in \{3,5,56\}$, so that $h_1(f_r(x))$ is irreducible. Observe that $x^2 - \frac{1}{m^2}$ maps the interval $(-1/m, 0)$ into itself, and in particular, $f_r^n(0) < 0$ for all $n \geq 1$. Thus $h_1(f_r^n(0)) < 0$ as well, and hence cannot be a square in $\Q$. Then Lemma \ref{squareirred} (with $g = h_1 \circ f_r$ and $f = f_r$) proves that $h_1(f_r^n(x))$ is irreducible for all $n \geq 2$. Thus for $s \not\in \{3,5,56\}$, we have that $h_1(f_r^n(x))$ is irreducible for all $n \geq 1$. We now present a result that builds on Corollary \ref{g2cor}.

\begin{corollary} \label{h2cor}
Let $r = 1/c$ and $c = -(s^2-1)^2$ for $s \geq 2$, and let $g_2$ be as in \eqref{gnotation} and $h_1,h_2$ as in \eqref{hdef}. Suppose that $(s^2-1)^2 \leq 10^9$. Then for all $n \geq 1$ we have $g_2(f_r^n(x))$ and $h_2(f_r^n(x))$ irreducible. If in addition $s \not\in \{3, 5, 56\}$ then for all $n \geq 1$ we have $h_1(f_r^n(x))$ irreducible. In particular if $(s^2-1)^2 \leq 10^9$ and $s \not\in \{3, 5, 56\}$, then $f_r^n(x)$ is a product of three irreducible factors for all $n \geq 2$.
\end{corollary}

\begin{proof}
Observe that $(s^2 - 1)^2 \leq 10^9$ if and only if $s \leq 177$. We have shown in Corollary \ref{g2cor} that $g_2(f_r^n(x))$ is irreducible for all $s$ with $2 \leq s \leq 177$. In the paragraph preceding the present corollary, we showed that $s \not\in \{3, 5, 56\}$ implies that $h_1(f_r^n(x))$ is irreducible for all $n \geq 1$. To show that $h_2(f_r^n(x))$ is irreducible for $n \geq 1$, it suffices by Lemma \ref{squareirred} to show that $\{-h_2(f_r(0))\} \cup \{h_2(f_r^n(0)) : n \geq 2\}$ contains no squares in $\Q$. Note that $-h_2(f_r(0)) =  -\frac{s^3 - s - 1}{(s^2-1)^2}$, and we have $s^3 - s - 1 > 0$ for $s \geq 2$. Hence $h_2(f_r(0))$ is not a square in $\Q$. To verify that $h_2(f_r^n(0))$ is a non-square in $\Q$ for all $n \geq 2$, we search for primes $p$ satisfying the condition \eqref{goal}, with $h_2$ replacing $g_2$. We find that there exists a prime $p \leq 500$ with the desired property for all $s$ with $2 \leq s \leq 177$ except for $s = 153$. For that $s$-value, the prime $p = 1051$ suffices.

For each such $s$ and $p$, we take the finitely many terms of the sequence $(h_2(f_r^n(0)))_{n \geq 2}$ that have still not been proven non-square, and reduce modulo small primes until all have been proven non-square.
Unsurprisingly, the $s$-value producing the largest number of such terms is $s = 153$, where we must check that each of $h_2(f_r(0)), h_2(f_r^2(0)), \ldots, h_2(f_r^{67}(0))$ is a non-square. In all cases the desired result is achieved by reducing modulo primes less than 100.
\end{proof}

Finally, we handle the case of $s \in \{3, 5, 56\}$. These are precisely the $s$-values for which $s^3 - s + 1$ is a square. In this case, $h_1(f(x))$ is no longer irreducible; indeed, we have
\begin{equation} \label{h12def}
h_1(f(x)) = \left(x - \frac{\sqrt{s^3 - s +1}}{s^2-1}\right)\left(x + \frac{\sqrt{s^3 - s +1}}{s^2-1}\right)=: h_{11}(x)h_{12}(x).
\end{equation}

\begin{proposition} \label{specials}
Let $r = 1/c$ and $c = -(s^2-1)^2$ for
$s \in \{3, 5, 56\}$. Let $g_2$ be as in \eqref{gnotation}, $h_2$ as in \eqref{hdef}, and $h_{11}$ and $h_{12}$ as in \eqref{h12def}. Then for all $n \geq 1$ we have $g_2(f_r^n(x))$, $h_2(f_r^n(x))$, $h_{11}(f_r^n(x))$, and $h_{12}(f_r^n(x))$ irreducible; in particular, $f_r^n(x)$ is a product of four irreducible factors for all $n \geq 3$.
\end{proposition}

\begin{proof}
Corollary \ref{h2cor} shows that for $s \in \{3, 5, 56\}$, we have $g_2(f_r^n(x))$ and $h_2(f_r^n(x))$ irreducible for all $n \geq 1$. To show that $h_{11}(f_r^n(x))$ and $h_{12}(f_r^n(x))$ are irreducible for $n \geq 1$, it suffices by Lemma \ref{squareirred} to show that none of
$$
\{-h_{11}(f_r(0))\} \cup \{-h_{12}(f_r(0))\} \cup \{ h_{11}(f_r^n(0)) : n \geq 2\} \cup \{h_{12}(f_r^n(0)) : n \geq 2\}
$$
is a square in $\Q$.
Note that $-h_{11}(f_r(0)) = ((s^2-1)(\sqrt{s^3 - s + 1}) + 1)/(s^2-1)^2$. For $s = 3, 5, 56$ respectively, the prime factorization of the numerator of $-h_{11}(f_r(0))$ is $41, 5 \cdot 53$, $2 \cdot 656783$, none of which is a square. Moreover, $-h_{12}(f_r(0)) < 0$, and hence cannot be a square. Also, one readily sees that $h_{11}(f_r^n(0)) < 0$ for all $n \geq 2$. For $s = 3$, we reduce the sequence $(h_{12}(f_r^n(0))_{n \geq 2}$ modulo 29 and find that it cycles among the four values $17, 15, 26, 21$, none of which is a square modulo $29$. For $s = 5$ we reduce modulo 23 and find that the sequence in question cycles between $10$ and $11$, which are both non-squares modulo $23$. For $s = 56$ we reduce modulo 31 and find that the sequence takes only the value 6, i.e. $h_{12}(f_r^n(0)) \equiv 6 \bmod{31}$ for all $n \geq 2$. But $6$ is non-square modulo 31.
\end{proof}

\subsection{The case of $f_r$ irreducible, $f_r^2$ reducible}
Assume now that $c = 4m^2(m^2-1)$ for some $m \geq 2$, in which case we have
\begin{equation} \label{qdef}
f_r^2(x) = \left(x^2-\frac{1}{m}x+\frac{2m^2-1}{4m^2(m^2-1)}\right)\left(x^2+\frac{1}{m}x+\frac{2m^2-1}{4m^2(m^2-1)}\right). 
\end{equation}
Let $q_1(x) = x^2-\frac{1}{m}x+\frac{2m^2-1}{4m^2(m^2-1)}$ and $q_2(x) = x^2+\frac{1}{m}x+\frac{2m^2-1}{4m^2(m^2-1)}.$ We note that $q_1$ and $q_2$ both have discriminant $-1/(m^2-1)$, and so are irreducible.

Observe that for $m = 2$ we have the factorization
\begin{equation} \label{m2fact}
q_2(f_r(x)) = (x^2 - (1/2)x + 19/48)(x^2 + (1/2)x + 19/48).
\end{equation}
However, this is the only $m$-value for which such a factorization occurs, as the next two results show.

\begin{proposition} \label{hyperell}
Let $r = 1/c$ and $c = 4m^2(m^2 - 1)$ for $m \geq 2$. If $f_3(x)$ has strictly more than two irreducible factors, then either
$$8m^6-12 m^4+4 m^3+4 m^2-4 m+1 \qquad \text{or} \qquad 8m^6-12 m^4-4 m^3+4 m^2+4 m+1$$
is a square in $\Q$.
\end{proposition}
\begin{proof} Observe that $f_r^3(x)$ has strictly more than two irreducible factors if and only if $q_i(f_r(x))$ is reducible for at least one $i \in \{1, 2\}$. Assume that $q_i(f_r(x))$ is reducible, let $\alpha$ be a root of $q_i(f_r(x))$, and observe that $f_r(\alpha) =: \beta$ is a root of $q_i(x)$. By the irreducibility of $q_i(x)$, we have $[\Q(\beta) : \Q] = 2$. Because
$q_i(f_r(x))$ is reducible, we have $[\Q(\alpha) : \Q] < 4$, which implies $[\Q(\alpha) : \Q(\beta)] = 1$, and thus $\alpha \in \Q(\beta)$.
But $\alpha$ is a root of $f_r(x) - \beta = x^2 + r -\beta$, and so $\alpha \in \Q(\beta)$ is equivalent to $\beta - r$ being a square in $\Q(\beta)$. Letting $\beta'$ be the other root of $q_i(x)$, we have
$$N_{\Q(\beta) / \Q} (\beta - r) = (\beta - r)(\beta' - r) = q_i(r) = \frac{8m^6 - 12m^4 \mp 4m^3 + 4m^2 \pm 4m + 1}{(4m^4 - 4m^2)^2}
$$
The multiplicativity of the norm map implies that the rightmost expression is a square in $\Q$.
\end{proof}

We now prove Theorem \ref{ellchab}, which we restate here.
\begin{theorem} \label{ellchab2}
The only rational points on the curve $y^2 = 8x^6 - 12x^4 - 4x^3 + 4x^2 + 4x + 1$ are those with $x \in \{-2, -1, 0, 1\}$.
\end{theorem}

\begin{proof} We note first that the map $(x,y)\rightarrow(1/x,y/x^3)$ gives a birational transformation from the curve $y^2=8x^6 - 12x^4 - 4x^3 + 4x^2 + 4x + 1$ to the curve 
\[C: y^2=F(x)=x^6 + 4x^5 + 4x^4 - 4x^3 - 12x^2 + 8.\] 
Therefore, it suffices to find all rational points on $C$. Next, we see that the polynomial $F(x)$ factors over a small extension of $\mathbb{Q}$. Fix an algebraic number $\beta$ satisfying $\beta^3 - 8\beta^2 + 20\beta - 8=0$, and observe that $F(x)$ factors as
\[\Big(x^2 + (-\beta + 4)x + 1/2(\beta^2 - 6\beta + 8)\Big)\Big(x^4 + \beta x^3 + 1/2(\beta^2 - 2\beta)x^2 - 4x - 2\beta + 4\Big).\]
Let $K = \mathbb{Q}(\beta)$, a number field of class number $1$. 
Therefore, if $(x,y)$ is a rational point on $C$, then there exist $y_1, y_2, \alpha \in K$ such that 
\begin{equation}\label{eq:split}
\begin{split} 
\alpha y_1^2&=F_1(x)=x^2 + (-\beta + 4)x + 1/2(\beta^2 - 6\beta + 8)\\[3pt] 
\alpha y_2^2&=F_2(x)=x^4 + \beta x^3 + 1/2(\beta^2 - 2\beta)x^2 - 4x - 2\beta + 4
\end{split} 
\end{equation} 
simultaneously; this follows from the fact that $F_1(x)$ and $F_2(x)$ lie in the same square-class in $K$. Moreover, we may assume that $\alpha$ is in the ring of integers $\mathcal{O}_K$ of $K$ and that the ideal $\alpha\mathcal{O}_K$ is not divisible by the square of an ideal in $\mathcal{O}_K$. On the other hand, since the degrees of $F_1$ and $F_2$ are not both odd (see Example 9 and Theorem 11 of \cite{stollrat}), if $\mathfrak{p}$ is a prime in $\mathcal{O}_K$ that divides $\alpha$ and is coprime to $2$, then $\mathfrak{p}$ must divide the resultant $R=36\beta^2 - 240\beta + 400$ of $F_1$ and $F_2$. Therefore, we may write  
\begin{equation}\label{possiblealphas} 
\alpha=(-1)^{e_0}\cdot2^{e_1}\cdot\Big(\frac{\beta^2}{4} - \frac{3\beta}{2} + 2\Big)^{e_2}\cdot\Big(\frac{3\beta^2}{4} - 4\beta + 5\Big)^{e_3}
\end{equation} 
for some $e_i\in\{0,1\}$ and $0\leq i\leq3$; here we use Sage to factor the fractional ideal generated by $R$ and find generators $-1$ and $\frac{\beta^2}{4} - \frac{3\beta}{2} + 2$ of the unit group of $K$. In particular, we have deduced that if $(x,y)\in C(\mathbb{Q})$, then $(x,y_2)$ is a $K$-point on 
\[V_\alpha: \alpha y^2=F_2(x),\] 
for some $y_2\in K$ and some $\alpha$ in (\ref{possiblealphas}). In particular, for such $\alpha$ it must be the case that $V_\alpha(K_v)$ is non-empty for every completion $K_v/K$. However, we check with MAGMA that only the curves $V_\alpha$ corresponding to $\alpha=1$ and $\alpha=\frac{\beta^2}{4} - \frac{3\beta}{2} + 2$ have points everywhere locally. On the other hand, $V_\alpha(K)$ is non-empty for both of these choices of $\alpha$. Therefore, there exist computable elliptic curves $E_1$ and $E_2$ (in Weierstrass form) together with birational maps $\phi_1: E_1\rightarrow V_1$ and $\phi_2: E_2\rightarrow V_{\frac{\beta^2}{4} - \frac{3\beta}{2} + 2}$ all defined over $K$. In particular, it suffices to compute the sets 
\[ S_i=\big\{P\in E_i(K): x(\phi_i(P))\in\mathbb{P}^1(\mathbb{Q})\big\}\] 
for $i\in\{1,2\}$, to classify the integral points on $C$. However, $E_1(K)$ and $E_2(K)$ both have rank $2$. In particular, $\text{rank}(E_1(K))$ and $\text{rank}(E_2(K))$ are both strictly less than $[K:\mathbb{Q}]=3$. Therefore, $S_1$ and $S_2$ are finite sets, and we may use the elliptic Chabauty method to describe them; see, for instance, \cite[\S4.2]{Bruin}. Moreover, since both $E_1$ and $E_2$ are in Weierstrass form and we succeed in finding explicit generators for their Mordell-Weil groups, we may use an implementation of the elliptic Chabauty method in MAGMA to describe $S_1$ and $S_2$; see the file named \emph{Elliptic Chabauty} at the website above for the relevant code. In particular, we deduce that 
\[C(\mathbb{Q})=\{\infty^+,\infty^-,(\pm{1},\pm{1}), (-1/2,\pm{19/8})\},\]
from which Theorem \ref{ellchab2} easily follows.     
\end{proof}

\begin{corollary} \label{chabcor}
Let $r = 1/c$ and $c = 4m^2(m^2-1)$ for $m \geq 2$. Then $f_r^3(x)$ has more than two irreducible factors if and only if $m = 2$.
\end{corollary}

\begin{proof}
The sufficiency is clear from \eqref{m2fact}. To see that $m = 2$ is also necessary, assume that $f_r^3(x)$ has more than two irreducible factors. From Proposition \ref{hyperell}, we have that $m$ or $-m$ is the $x$-coordinate of an integral point on the curve $y^2 = 8x^6 - 12x^4 - 4x^3 + 4x^2 + 4x + 1$. It then follows from Theorem \ref{ellchab2} that $\pm m \in \{-2, -1, 0, 1\}$. Since $m \geq 2$, the only possibility is $m = 2$.
\end{proof}

We have now assembled enough ingredients to prove Theorem \ref{maincases}.

\begin{proof}[Proof of Theorem \ref{maincases}]
Part (a) is proven in Propositions \ref{m4prop} and \ref{m4}. Part (b) follows from Proposition \ref{specials} and the remarks after \eqref{hdef}.

 The first assertion of part (c) is proven in Corollary \ref{chabcor}. To prove the second assertion, let $m = 2$, let $q_1$ and $q_2$ be as in \eqref{qdef}, and set
$v_1(x) = x^2 - (1/2)x + 19/48$ and $v_2(x) = x^2 + (1/2)x + 19/48$, so that $q_2(f_r(x)) = v_1(x)v_2(x)$. We must show that $q_1(f_r^n(x))$ and $v_j(f_r^n(x))$ ($j \in \{1, 2\}$) are irreducible for all $n \geq 1$. Because $q_1, v_1$, and $v_2$ have even degree, by Lemma \ref{squareirred} it suffices to prove $q_1(f_r^n(0))$ and $v_j(f_r^n(0))$ are not squares in $\Q$ for all $n \geq 1$.

We now search for primes $p$ satisfying the condition \eqref{goal}, with $q_1$ and $v_j$ replacing $g_2$. We reduce the sequence $q_1(f_r ^n(0))$ modulo $239$, and find that it only takes the non-square value $13$ for $n \geq 7$. For $n$ with $1 \leq n \leq 6$, one verifies directly that $q_1(f_r ^n(0))$ is not a square.
We reduce the sequence $v_1(f_r ^n(0))$ modulo $239$, and find that it only takes the non-square value $73$ for $n \geq 7$. For $n$ with $1 \leq n \leq 6$, one verifies directly that $v_1(f_r ^n(0))$ is not a square. We reduce the sequence $v_2(f_r ^n(0))$ modulo $41$, and find that it only takes the non-square value $24$ for $n \geq 7$. For $n$ with $1 \leq n \leq 6$, one verifies directly that $v_2(f_r ^n(0))$ is not a square.
\end{proof}

We close this section with a proof of one case of Theorem \ref{bigcomp}.

\begin{proposition} \label{h22cor}
Let $r = 1/c$ and $c = 4m^2(m^2-1)$ for $m \geq 3$, and let $q_1$ and $q_2$ be as in \eqref{qdef}. Suppose that $4m^2(m^2-1) \leq 10^9$. Then for all $n \geq 1$ we have $q_1(f_r^n(x))$ and $q_2(f_r^n(x))$ irreducible. Hence $f_r^n(x)$ is a product of two irreducible factors for all $n \geq 2$.
\end{proposition}

\begin{proof}
Observe that $4m^2(m^2-1) \leq 10^9$ if and only if $m \leq 125$. Because $q_1$ and $q_2$ have even degree, by Lemma \ref{squareirred} it suffices to prove $q_1(f_r^n(0))$ and $q_2(f_r^n(0))$ are non-squares in $\Q$ for all $n \geq 1$. We search for primes $p$ satisfying the condition \eqref{goal}, with $q_1$ and $q_2$ replacing $g_2$.

For $q_1(f_r^n(0))$, we find that there exists a prime $p \leq 500$ (indeed, $p \leq 337$) with the desired property for all $m$ with $3 \leq m \leq 125$. For $q_2(f_r^n(0))$, we also find that there exists a prime $p \leq 500$ with the desired property for all $m$ with $3 \leq m \leq 125$.

For each such $m$ and $p$, we take the finitely many terms of the sequence $(q_1(f_r^n(0)))_{n \geq 2}$ (resp. $(q_2(f_r^n(0)))_{n \geq 2}$) that have still not been proven non-square, and reduce modulo small primes until all have been proven non-square.
\end{proof}

\section{The proof of cases (1)-(4) of Theorem \ref{main}} \label{preliminaries}

In the last section, we saw the primary importance of whether or not $p(f_r^n(0))$ is a square, for various polynomials $p(x)$. For the remainder of this article, we use similar ideas to study the irreducibility of $f_r(x)$ in the case where $f_r^2(x)$ is irreducible. However, we use a refinement of \cite[Proposition 4.2]{quaddiv}, similar to \cite[Theorem 2.3]{itconst}, that is more powerful; see Lemma \ref{biglemma} (restated as Lemma \ref{biglemma2} below).

Recall from the introduction that $r = 1/c$, and that $f_r^n(0)$ is a rational number with denominator $c^{2^{n-1}}$. We define $a_n(c)$ to be the numerator of $f_r^n(0)$. Hence $a_n(c)$ is described by the recurrence
\begin{equation} \label{andef2}
a_1(c) = 1, \qquad a_{n}(c) = a_{n-1}(c)^2 + c^{2^{n-1} - 1} \quad \text{for $n \geq 2$}.
\end{equation}
To ease notation, we often suppress the dependence on $c$, and write $a_1, a_2,$ etc. Recall also that we define
\begin{equation} \label{bndef2}
b_n := \frac{a_{n-1} + \sqrt{a_{n}}}{2} \in \overline{\mathbb{Q}}.
\end{equation}

\begin{proposition} \label{negative}
If $c < 0$, then $a_n$ is not a square in $\Q$ for all $n \geq 2$.
\end{proposition}

\begin{proof}
Let $r = 1/c$ and $f_r(x) = x^2 + r$, and consider the image of the interval $I=(-\sqrt{-r},0)$ under $f_{r}:\R\to \R$. We have $f_{r}(-\sqrt{-r})=0$ and $f_{r}(0)=r\in I$, so as $f_{r}$ is a continuous function with no critical points in $I$, it follows that $f_{r}(I)\subset I$. As $f_{r}(0)=r\in I$,  inductively, $f^n_{r}(0)\in I$ for all $n\geq 1$. Hence $0 > f_r^n(0) = a_n/c^{2^{n}}$, and hence $a_n < 0$ for $n \geq 1$, proving that $a_n$ is not a square in $\Q$.
\end{proof}

We now prove Lemma \ref{biglemma}, which we restate here.
\begin{lemma} \label{biglemma2}
Suppose that $c \in \Z \setminus \{0\},$ $r = 1/c$, and $f_r^2$ is irreducible. Let $a_n = a_n(c)$ and $b_n$ be defined as in \eqref{andef2} and \eqref{bndef2}, respectively.
If for every $n \geq 3$, $b_n$ is not a square in $\Q$ (which holds in particular if $a_n$ is not a square in $\Q$), then $f_r^n(x)$ is irreducible for all $n \geq 1$.
\end{lemma}

\begin{proof} This proof is essentially the same as the proof of \cite[Theorem 2.3]{itconst}, but for completeness we give the argument here.
By hypothesis $f_r^2(x)$ is irreducible; assume inductively that $f^n_r(x)$ is irreducible for some $n \geq 2$. Let $\alpha$ be a root of $f_r^{n+1}(x)$, and observe that $f_r(\alpha) =: \beta$ is a root of $f_r^n(x)$. By our inductive assumption, we have $[\Q(\beta) : \Q] = 2^n$. Now $f_r^{n+1}(x)$ is irreducible if and only if $[\Q(\alpha) : \Q] = 2^{n+1}$, which is equivalent to $[\Q(\alpha) : \Q(\beta)] = 2$. This holds if and only if $f_r(x) - \beta$ is irreducible over $\Q(\beta)$, i.e. $\beta - r$ is not a square in $\Q(\beta)$. Now factor
$f_r^n(x)$ over $K_1 := \Q(\sqrt{-r})$. We have $f_r^n(x) = (f_r^{n-1}(x) - \sqrt{-r})(f_r^{n-1}(x) + \sqrt{-r})$, and because $[\Q(\beta) : \Q] = 2^n$, we must have $[\Q(\beta) : K_1] = 2^{n-1}$, which implies that the minimal polynomial of $\beta$ over $K_1$ is one of $f_r^{n-1}(x) \pm \sqrt{-r}$. It follows that $N_{\Q(\beta)/K_1}(\beta - r)$ is the product of $(\beta' - r)$, where $\beta'$ varies over all roots of $f_r^{n-1}(x) \pm \sqrt{-r}$; this product is just $f_r^{n-1}(r) \pm \sqrt{-r}$ (here we use that $n \geq 2$, so the degree of $f_r^{n-1}(x)$ is even and we may replace the product of $(\beta' - r)$ with the product of $(r-\beta')$). To summarize, we have
$$
N_{\Q(\beta)/K_1}(\beta - r) = f_r^{n-1}(r) \pm \sqrt{-r} = f_r^{n}(0) \pm \sqrt{-r}.
$$
Suppose now that $f_r^{n+1}(x)$ is reducible, and hence $\beta-r$ is a square in $\Q(\beta)$. Because the norm map is multiplicative, this implies $N_{\Q(\beta)/K_1}(\beta - r)$ is a square in $K_1$, i.e. there exist $s_1, s_2 \in \Q$ with
$(s_1 + s_2 \sqrt{-r})^2 = f_r^{n}(0) \pm \sqrt{-r}$. Elementary calculations show this last equality implies $s_2=\frac{1}{2s_1}$ and $s_1^2 - rs_2^2 = f_r^n(0)$, whence $$s_1^2=\frac{f_r^{n}(0)\pm\sqrt{f_r^{n+1}(0)}}{2} = \frac{a_{n}\pm \sqrt{a_{n+1}}}{2c^{2^{n-1}}}.$$
Now $n \geq 2$, and hence we have that one of $(a_{n}\pm \sqrt{a_{n+1}})/2$ is a square in $\Q$. If $c < 0$, then this is impossible by Proposition \ref{negative}. Hence suppose $c > 0$.
As $a_{n+1}=a_{n}^2+c^{2^{n}-1}>a_{n}^2>0$, we have $(a_{n}-\sqrt{a_{n+1}})/2<0$, implying that $(a_{n} + \sqrt{a_{n+1}})/2$ is a square in $\Q$. But this is contrary to the hypotheses of the lemma, and we thus conclude that $f_r^{n+1}(x)$ is irreducible.
\end{proof}

\begin{proposition}
\label{a3nonsquare}
Let $c \in \Z \setminus \{0, -1\}$. Then neither $a_3$ nor $a_4$ is a square in $\Q$.
\end{proposition}
\begin{proof}
We have $a_3(c) = c^3 + c^2 + 2c + 1$, and so if $a_3(c) = y_0^2$ for $y_0 \in \Q$, then necessarily $y_0 \in \Z$, and $(c, y_0)$ is an integer point on the elliptic curve $y^2 = x^3 + x^2 + 2x + 1$. This curve has conductor 92, and is curve 92.b2 in the LMFDB \cite{lmfdb}. Besides the point at infinity, it has only the rational points $(0, \pm 1)$, but $c = 0$ is excluded by hypothesis.

We now address $a_4(c)$. As in the previous paragraph, if $a_4(c) = y_0^2$ for $y_0 \in \Q$, then $(c,y_0)$ is an integer point on the hyperelliptic curve
\begin{equation*}
C : y^2 = x^7 + x^6 + 2x^5 + 5x^4 + 6x^3 + 6x^2 + 4x + 1.
\end{equation*}
One easily checks that $x^7 + x^6 + 2x^5 + 5x^4 + 6x^3 + 6x^2 + 4x + 1$ has no repeated roots, and hence $C$ has genus 3. Denote by $J$ the Jacobian of $C$. A two-descent using MAGMA \cite{magma} shows that $J(\Q)$ has rank zero, and hence consists only of torsion. We now use standard reduction techniques to determine all torsion in $J(\Q)$ \cite[Theorem C.1.4 and Section C.2]{jhsdioph}. We have a commutative diagram
\begin{equation} \label{commdiag}
\begin{CD}
C(\Q) @>>> J(\Q)\\
@VVV @VVV \\
C(\mathbb{F}_3) @>>> J(\mathbb{F}_3)
\end{CD}
\end{equation}
where the vertical maps are reduction modulo 3 and the horizontal maps are the Abel-Jacobi maps taking $P$ to the divisor class of $(P - \infty)$. The latter are injective \cite[Corollary A.6.3.3]{jhsdioph}. The discriminant of $C$ is $2^{12} \cdot 23 \cdot 2551$, and it follows that $C$, and hence $J$ \cite[p. 164]{jhsdioph}, has good reduction at all primes $p \notin \{2, 23, 2551\}$. Because $J(\Q)$ is torsion, it follows that for any such prime $p$, the reduction map $J(\Q) \to J(\mathbb{F}_{p})$ is injective; see, for instance, the appendix of \cite{katz}. Thus the right vertical map in \eqref{commdiag} is injective, and it follows that the left vertical map is injective as well. But one verifies that $\#C(\mathbb{F}_3) = 4$, and hence
$C(\Q) = \{\infty, (0,\pm 1), (-1,0)\}$. Because we have excluded $c = 0, -1$, we arrive at the desired contradiction.

One may attempt the same argument with $a_5(c)$, but a 2-descent on the Jacobian $J$ of the associated genus-7 hyperelliptic curve shows only that the rank of $J(\Q)$ is at most 2.
\end{proof}

\begin{proposition} \label{divisibility}
The sequence $(a_n)_{n \geq 1}$ is a rigid divisibility sequence. (See Definition \ref{rigid}).
\end{proposition}

\begin{proof}
This is a straightforward application of \cite[Lemma 2.5]{zdc}.
\end{proof}

\begin{proposition} \label{congruences}
Suppose that $c+1$ is not a square in $\Z$. If $c$ satisfies any of the congruences in Table 1, then $a_n$ is not a square in $\Q$ for all $n \geq 2$.
\end{proposition}
\begin{table} \label{congtable}
\begin{align*}
c & \equiv 1,2 & \pmod{3} \\
c & \equiv 3 & \pmod{4} \\
c & \equiv 2,3 & \pmod{5} \\
c & \equiv 1,2,5,6 & \pmod{7}\\
c & \equiv 1 & \pmod{8}\\
c & \equiv 1,3,5,7,10 & \pmod{11} \\
c & \equiv 3, 4, 5, 6, 8, 11 & \pmod{13}\\
c & \equiv 6, 10, 14, 15 & \pmod{17} \\
c & \equiv 1, 4, 9, 11, 12, 13, 15, 16, 18 & \pmod{19}\\
c & \equiv 6, 10, 12, 18, 20, 22 & \pmod{23} \\
c & \equiv 2, 12, 14, 17, 18, 27 & \pmod{29}\\
c & \equiv 1, 10, 13, 16, 22, 27, 30 & \pmod{31} \\
c & \equiv 6, 18, 23, 31, 32, 35 & \pmod{37}\\
c & \equiv 7, 8, 11, 19, 25, 28, 35, 36 & \pmod{41} \\
c & \equiv 1, 2, 4, 5, 9, 14, 15, 21, 27, 33, 37, 42 & \pmod{43}\\
c & \equiv 6, 7, 9, 10, 24, 25, 28, 33, 46 & \pmod{47} \\
c & \equiv 5, 18, 21, 23, 26, 30, 37, 40, 43, 45, 46, 47 & \pmod{53}\\
c & \equiv 10, 14, 16, 29, 37, 47, 55, 57, 58 & \pmod{59} \\
c & \equiv 2, 3, 11, 13, 15, 27, 30, 32, 34, 40, 45, 50 & \pmod{61}\\
c & \equiv 10, 15, 20, 32, 33, 38, 41, 49, 51, 53, 55, 66 & \pmod{67} \\
c & \equiv 4, 10, 49, 51, 53, 61, 70 & \pmod{71} \\
c & \equiv 1, 3, 35, 43, 44, 50, 51, 71 & \pmod{73} \\
c & \equiv 3, 12, 25, 32, 36, 58, 78 & \pmod{79} \\
c & \equiv 15, 16, 19, 23, 25, 29, 31, 37, 41, 44, 51, 56, 59, 68, 71, 82 & \pmod{83} \\
c & \equiv 13, 25, 49, 63 & \pmod{89} \\
c & \equiv 3, 9, 21, 53, 59, 79, 89 & \pmod{97} \\
\end{align*}
\caption{Congruences that ensure $a_n$ is not a square for $n \geq 2$, provided that $c+1$ is not a square.}
\end{table}

\begin{proof}
By Proposition \ref{negative}, it suffices to consider $c > 0$. Because $a_2 = c+1 > 0$ is non-square by assumption, there is a prime $q$ with $v_q(c+1)$ odd. Proposition \ref{divisibility} then implies that $a_{2m}$ is non-square for all $m\geq 1$, so we need only check that $a_n$ is non-square for odd $n \geq 2$. To do this, we let $f(x) = x^2 + 1/c$ and we take $p$ to be a fixed prime with $p < 100$ and $p \nmid c$. Let $c_0 \in \{1, \ldots, p-1\}$ satisfy $(1/c) \equiv c_0 \bmod{p}$ and put $\bar{f} = x^2 + c_0 \in \mathbb{F}_p[x]$. Now $a_n = c^{2^{n-1}}f^n(0)$, and it follows that if $\bar{f}^n(0)$ is not a square in $\mathbb{F}_p$, then $a_n$ is not a square in $\Q$. The sequence $(\bar{f}^n(0) \bmod{p})_{n \geq 1}$ eventually lands in a repeating cycle. When this sequence is such that $\bar{f}^{2n+1}(0)$ is a non-square in $\mathbb{F}_p$ for all $n \geq 2$, then $a_{2n+1}$ is a non-square in $\Z$ for all $n \geq 1$ (the $n=1$ case is by Proposition \ref{a3nonsquare}). Most of the pairs of $p,c$ listed in Table 1 yield such a result. For instance, when $p = 3$ and $c \equiv 1 \bmod{p}$, we have
$\bar{f}^{n}(0) = 2$ for all $n \geq 2$. When $p = 5$ and $c \equiv 3 \bmod{p}$, the sequence $\bar{f}^n(0)$ is $2, 1, 3, 1, 3, \ldots$, and hence $\bar{f}^{2n+1}(0)$ is a non-square for all $n \geq 1$. The remaining pairs $p,c$ in Table 1 satisfy the condition that both $\bar{f}^{3(n+1)+1}(0)$ and $\bar{f}^{3n+2}(0)$ are non-squares for $n \geq 1$ (the $n+1$ comes from the fact that $a_4$ is automatically a non-square by Proposition \ref{a3nonsquare}). Thus $a_{3n+1}$ and $a_{3n+2}$ are non-squares in $\Z$ for all $n \geq 1$. But by Proposition \ref{a3nonsquare} we have that $a_3$ is not a square in $\Z$, and it follows from Proposition \ref{divisibility} that $a_{3n}$ is a non-square in $\Z$ for all $n \geq 1$. An example is when $p = 7$ and $c \equiv 5 \bmod{p}$, for which the sequence $\bar{f}^n(0)$ is $3, 5, 0, 3, 5, 0, \ldots$.
\end{proof}

We now prove cases (1)-(4) of Theorem \ref{main}, which we restate here.
\begin{theorem} \label{main13}
Let $f_r(x) = x^2 + r$ with $r = 1/c$ for $c \in \Z \setminus \{0, -1\}$, and let $a_n$ and $b_n$ be as in \eqref{andef2} and \eqref{bndef2}. Assume that $c$ satisfies one of the following conditions:
\begin{enumerate}
\item $-c \in \Z \setminus \Z^2$ and $c < 0$;
\item $-c, c+1 \in \Z \setminus \Z^2$ and $c \equiv -1 \bmod{p}$ for a prime $p \equiv 3 \bmod{4}$;
\item $-c, c+1 \in \Z \setminus \Z^2$ and $c$ satisfies one of the congruences in Proposition \ref{congruences} (see Table 1);
\item $-c \in \Z \setminus \Z^2$ and $c$ is odd;
\end{enumerate}
In cases (1)-(3), $a_n$ is not a square in $\Q$ for any $n \geq 2$, while in case (4), $b_n$ is not a square for any $n \geq 2$. In all cases, $f_r^n(x)$ is irreducible for all $n \geq 1$.
\end{theorem}
\begin{proof}
Observe that conditions (1)-(4) each imply that $f_r^2(x)$ is irreducible, by Proposition \ref{irredprop} (note that $c = 4m^2(m^2 - 1)$ implies that $c+1 = (2m-1)^2$, and that this is impossible when $c$ is odd). We now argue that in cases (1)-(3) $a_n$ is not a square in $\Q$ for any $n \geq 2$ and in case (4), $b_n$ is not a square for any $n \geq 2$. In all these cases, Lemma \ref{biglemma} proves that $f_r^n(x)$ is irreducible for all $n \geq 1$.

If we are in case (1), then the desired conclusion holds by Proposition \ref{negative}.

Assume we are in case (2). Because we have already established case (1), it suffices to consider $c > 0$. Because $1/c \equiv -1 \bmod{p}$, we see that modulo $p$, the orbit of $0$ under $f_r$ is $0 \mapsto -1 \mapsto 0 \mapsto \cdots$. Moreover, $-1$ is not a square modulo $p$ by assumption, and so $a_{2n+1}$ is not a square for all $n \geq 3$. Because $a_2 = c+1 \geq 2$ is assumed non-square, it must be divisible by some prime to odd multiplicity. From Proposition \ref{divisibility} it then follows that $a_{2n}$ is not a square in $\Q$ for all $n \geq 1$.

In case (3) the desired conclusion holds by Proposition \ref{congruences}.

In case (4), if $a_n$ is not a square in $\Q$ then $b_n$ cannot be a square in $\Q$, and so we are done. If $a_n$ is square in $\Q$, then from the recursion in \eqref{andef2} and the fact that any integer equals its square modulo 2, we have
$$
\sqrt{a_n} \equiv a_n \equiv a_{n-1}^2+c^{2^{n-1} - 1} \equiv a_{n-1}^2+1 \equiv a_{n-1}+1 \pmod{2}.
$$
Thus modulo $2$, we have $a_{n-1}+\sqrt{a_n}\equiv 2a_{n-1}+1\equiv 1$, whence $v_2\left(\frac{a_{n-1}+\sqrt{a_n}}{2}\right)=-1$, proving that $b_n = \frac{a_{n-1}+\sqrt{a_n}}{2}$ is not a square in $\Q$.
\end{proof}

\section{Proof of cases (5) and (6) of Theorem \ref{main}} \label{mainresults}

In this section we deduce consequences from the assumption that $a_n(c)$
or even $b_n(c)$ is a square. This will lead to a fairly small upper bound
on~$n$ in terms of~$c$. One application is the proof of cases (5) and~(6)
of Theorem~\ref{main}. Another is the development of a fast algorithm
for checking that all iterates of~$f$ are irreducible as soon as $f^2$ is,
for all~$c$ up to a very large bound; this is done in the next section.

We denote the set of positive integers by $\Zplus$.

\begin{lemma} \label{L:dioph}
  Let $c \in \Zplus$ and $n \ge 2$ such that $a_n(c)$ is a square.
  Then we can write $c = u v$ with coprime integers $u$ and~$v$ such that
  \begin{enumerate}
    \item if $c$ is odd, then
          \[ v^{2^{n-1}-1} - u^{2^{n-1}-1} = 2 a_{n-1}(uv) ; \]
    \item if $c$ is even, then $u$ is even and
          \[ v^{2^{n-1}-1} - \tfrac{1}{4} u^{2^{n-1}-1} = a_{n-1}(uv) . \]
  \end{enumerate}
  If in addition $b_n(c) = (a_{n-1}(c) + \sqrt{a_n(c)})/2$ is a square (with the positive
  square root), then $c$ is even and $v$ is a square (and $u$ and~$v$ are positive)
  or $-u$ is a square (and $u$ and~$v$ are negative).
\end{lemma}

\begin{proof}
  To simplify notation, we set $N := 2^{n-1} - 1$.
  By assumption, there is $s \in \Zplus$ such that
  \[ a_n(c) = c^N + a_{n-1}(c)^2 = s^2 \]
  and hence
  \[ c^N = \bigl(s + a_{n-1}(c)\bigr) \bigl(s - a_{n-1}(c)\bigr) . \]
  It follows easily by induction that $a_m(c) \equiv 1 \bmod c$ for all~$m \ge 1$;
  in particular, $a_{n-1}(c)$ and~$s$ are coprime with~$c$. Since
  $\gcd(a_{n-1}(c), s)$ divides a power of~$c$, it follows that $a_{n-1}(c)$
  and~$s$ are also coprime. So we can deduce that
  \[ \gcd\bigl(s + a_{n-1}(c), s - a_{n-1}(c)\bigr) \mid \gcd(2s, 2a_{n-1}(c)) = 2 . \]
  We set $t_+ := s + a_{n-1}(c)$ and $t_- := s - a_{n-1}(c)$.
  \begin{enumerate}
    \item If $c$ is odd, then the gcd on the left is odd (since it divides a power of~$c$),
          so $t_+$ and~$t_-$ are coprime.
          Then $t_+ t_- = c^N$ implies that $c = uv$ with $u$, $v$ coprime
          and $t_+ = v^N$, $t_- = u^N$. The claim follows, since $t_+ - t_- = 2 a_{n-1}(c)$.
    \item Now assume that $c$ is even. Then $\gcd(s + a_{n-1}(c), s - a_{n-1}(c)) = 2$,
          since both entries have the same parity and their product is even.
          We can then write $c = uv$ with coprime $u$ and~$v$ and $u$~even such that
          either $t_+ = 2v^N$ and $t_- = \tfrac{1}{2}u^N$ or $t_+ = \tfrac{1}{2}u^N$ and $t_- = 2v^N$.
          In the first case, the claim again follows from $t_+ - t_- = 2 a_{n-1}(c)$.
          In the second case, we obtain
          $(-v)^N - \tfrac{1}{4} (-u)^N = (-t_- + t_+)/2 = a_{n-1}((-u)(-v))$,
          so we get the claim upon changing the signs of $u$ and~$v$.
  \end{enumerate}
   For the last claim, observe that
  \[ 0 < \frac{a_{n-1}(c) + \sqrt{a_n(c)}}{2} = \frac{a_{n-1}(c) + s}{2} = \frac{t_+}{2} . \]
  If $c$ is odd, then $t_+$ is odd, and $t_+/2$ cannot be a square. Otherwise, $t_+/2$
  is equal to either $v^N$ or $(-u)^N/4$. Since $N$ is odd, the claim follows.
\end{proof}

We set, for $c \ge 4$,
\[ F(c) = \frac{1}{2}\Bigl(1 - \sqrt{1 - \frac{4}{c}}\Bigr)
        = \frac{2}{c}\Bigl(1 + \sqrt{1 - \frac{4}{c}}\Bigr)^{-1} .
\]
From the first expression, it is clear that $F(c)$ decreases monotonically
from~$1/2$ to~$0$ as $c$ grows from~$4$ to infinity. The second expression
shows that for large~$c$, $F(c)$ is close to~$1/c$.

\begin{lemma} \label{L:anasymp}
  Let $c \ge 4$. Then the sequence $(\bar{a}_n(c))_{n \ge 1}$, where
  \[ \bar{a}_n(c) = \frac{a_n(c)}{c^{2^{n-1}-1}}, \]
  satisfies $1 = \bar{a}_1(c) < \bar{a}_2(c) < \ldots {}$ and
  $\lim_{n \to \infty} \bar{a}_n(c) = c F(c)$.
\end{lemma}

\begin{proof}
  We have that $\bar{a}_{n+1}(c) = 1 + \bar{a}_n(c)^2/c$. When $1 \le x < c F(c)$,
  then $c F(c) > 1 + x^2/c > x$, so that the sequence is strictly increasing and bounded
  by~$c F(c)$. Since $c F(c)$ is the smallest fixed point~$\ge 1$ of $x \mapsto 1 + x^2/c$,
  it must be the limit.
\end{proof}

We make a couple of definitions.

\begin{definition}
  Let $c \ge 2$ be an integer. We set
  \[ q(c) = \min \Bigl\{ \frac{v}{u} : \text{$u, v \in \Zplus$ coprime with $v > u$ and $c = uv$} \Bigr\} \]
  and
  \begin{multline} \tilde{q}(c) = \min \Bigl\{ \frac{v}{u} : \text{$u, v \in \Zplus$ coprime with $v > u$, $c = uv$,} \\
                                                   \text{and at least one of $u$ and $v$ is a square} \Bigr\} .
  \end{multline}
\end{definition}

We note that $\tilde{q}(c) \ge q(c) > 1 + 1/\sqrt{c}$, since $v \ge u + 1$ in the set above, so
$q(c) \ge 1 + 1/u$ for the minimizing~$u$, and $u < \sqrt{c}$, since $u^2 < uv = c$.

We write ``$\log$'' for the natural logarithm.

\begin{definition}
  Let $c \in \Zplus$ and $n \ge 2$. We define $\eps(n, c)$ so that
  \[ \log \frac{\sqrt{a_n(c)} + a_{n-1}(c)}{\sqrt{a_n(c)} - a_{n-1}(c)} = \frac{\eps(n, c)}{\sqrt{c}} . \]
\end{definition}

It follows from Lemma~\ref{L:anasymp} and the properties of~$F(c)$
that for fixed $c \ge 4$, $\eps(n, c)$ increases with~$n$ with limit
\[ \eps(c) := \lim_{n \to \infty} \eps(n, c) = \sqrt{c} \log \frac{1 + \sqrt{F(c)}}{1 - \sqrt{F(c)}} \] \label{epsdef}
and that $\eps(c)$ decreases monotonically when $c$ increases,
with $\lim_{c \to \infty} \eps(c) = 2$.
 In particular, we have that
\[ \eps(n, c) \le \eps(c) \le \eps(4) = 4 \log (1 + \sqrt{2}) \qquad\text{and}\qquad
   \frac{\eps(n,c)}{\sqrt{c}} \le 2 \log (1 + \sqrt{2}) \,.
\]
Since $(e^x-1)/x$ is monotonically increasing for positive~$x$, this implies that (for $c \ge 4$)
\begin{equation} \label{E:ineq}
  \exp\Bigl(\frac{\eps(n,c)}{\sqrt{c}}\Bigr) \le 1 + \frac{1+\sqrt{2}}{\log(1+\sqrt{2})} \cdot \frac{\eps(n,c)}{\sqrt{c}}
                                             \le 1 + \frac{4(1+\sqrt{2})}{\sqrt{c}} \,.
\end{equation}

We note that
\begin{align}
  \frac{\eps(c)}{\sqrt{c} \log q (c)} &< 3.46 & &\text{for $c \ge 4$,} \label{E:c4} \\
  \frac{\eps(c)}{\sqrt{c} \log(1 + 1/\sqrt{c})} &< 2.12 & &\text{for $c \ge 100$,} \label{E:c100} \\
  \frac{\eps(c)}{\sqrt{c} \log(1 + 1/\sqrt{c})} &< 2.01 & &\text{for $c \ge 10400$.} \label{E:c10400}
\end{align}
(To get~\eqref{E:c4}, we use \eqref{E:c100} and the explicit values of~$q(c)$
for $c < 100$. The maximum is achieved for $c = 6$.)
We will also need the elementary bound
\begin{equation} \label{E:logbound}
  \frac{1}{\log(1 + 1/\sqrt{c})} \le \sqrt{c} + \tfrac{1}{2} \,.
\end{equation}

We can now deduce an upper bound on~$n$ such that $a_n(c)$ can be a square.

\begin{proposition} \label{Prop}
  Let $c \ge 4$ be an integer and $n \ge 4$. If $c$ is odd or
  \[ n \ge 1 + \log_2\Bigl(1 + \frac{\eps(n,c)}{\sqrt{c} \log q(c)} + \frac{\log 4}{\log q(c)}\Bigr) , \]
  then $a_n(c)$ is not a square. This is the case whenever
  \[ \sqrt{c} \le \frac{2^{n-1}-1}{\log 4} - 3 . \]
  If the weaker condition
  \[ n \ge 1 + \log_2\Bigl(1 + \frac{\eps(n,c)}{\sqrt{c} \log \tilde{q}(c)}
                             + \frac{\log 4}{\log \tilde{q}(c)}\Bigr)
  \]
  holds, then $b_n(c)$ is not a square.
\end{proposition}

\begin{proof}
  In the following, we write $a_m$ for~$a_m(c)$, since $c$ is  fixed.
  We assume that $a_n$ is a square, so by  Proposition~\ref{a3nonsquare}
  we have $n \ge 5$, and by  Lemma~\ref{L:dioph} and its proof we can
  write $c = uv$ with coprime $u$, $v$ (and $u$ even when $c$ is even)
  such that
  \begin{align*}
    (v^N, u^N) &= (\sqrt{a_n} + a_{n-1}, \sqrt{a_n} - a_{n-1}) & & \text{if $c$ is odd;} \\
    (v^N, u^N) &= \bigl(\tfrac{1}{2}(\pm \sqrt{a_n} + a_{n-1}), 2(\pm\sqrt{a_n} - a_{n-1})\bigr)
                 & & \text{if $c$ is even,}
  \end{align*}
  where $N = 2^{n-1}-1$.

  First assume that $c$ is odd. Then $v > u > 0$, and we obtain  using~\eqref{E:ineq}
  \[ 1 + \frac{N}{\sqrt{c}} \le \Bigl(1 + \frac{1}{\sqrt{c}}\Bigr)^N
        < \Bigl(\frac{v}{u}\Bigr)^N
        = \frac{\sqrt{a_n} + a_{n-1}}{\sqrt{a_n} - a_{n-1}}
        = \exp\Bigl(\frac{\eps(n,c)}{\sqrt{c}}\Bigr)
         \le 1 + \frac{4(1+\sqrt{2})}{\sqrt{c}} ,
  \]
  which  is a contradiction, since $N \ge 15$.
  So $a_n$ cannot be a square.

  Now assume that $c$ is even. If $u, v > 0$ (this corresponds to the
  positive sign above), then
  \[ \Bigl(\frac{v}{u}\Bigr)^N = \frac{1}{4} \frac{\sqrt{a_n} + a_{n-1}}{\sqrt{a_n} - a_{n-1}} . \]
  If $u, v < 0$, then
  \[ \Bigl(\frac{u}{v}\Bigr)^N = 4 \frac{\sqrt{a_n} + a_{n-1}}{\sqrt{a_n} - a_{n-1}} . \]
  In both cases, we have that $|\log (v/u)| \ge \log q(c)$. This gives
  \begin{equation} \label{E:uveps}
    N \log q(c) \le  N \Bigl|\log \frac{v}{u}\Bigr| \le \log 4 + \frac{\eps(n,c)}{\sqrt{c}} ,
  \end{equation}
  which is equivalent to the inequality we wanted to show. If we assume that
  $b_n(c)$ is a square, then we have in addition that $|u|$ or~$|v|$ is a square,
  hence the bound is valid for~$\tilde{q}(c)$ in place of~$q(c)$.

  The bound on~$\sqrt{c}$ follows from the first inequality, the estimates
  $\eps(n, c) \le \eps(c)$, $q(c) > 1 + 1/\sqrt{c}$, and from \eqref{E:c4}
  and~\eqref{E:logbound}. Note that $3.46/(\log 4) + 0.5 < 3$.
\end{proof}

This gives the following.

\begin{corollary} \label{C:improve5.4}
  Let $c \ge 4$ be an integer and set $f(x) = x^2 + 1/c$.
  \begin{enumerate}
    \item If $c$ is odd, then all $f^n$ are irreducible.
    \item If $c$ is even and $f^m$ is irreducible for
          \[ m = 1 + \Bigl\lfloor \log_2\Bigl(1 + \frac{\log 4 + \eps(c)/\sqrt{c}}%
                                                       {\log(1 + 1/\sqrt{c})}\Bigr) \Bigr\rfloor ,
          \]
          then all $f^n$ are irreducible.
    \item If $c$ is even, $f^2$ is irreducible, and $a_p(c)$ is not a square
          for all prime numbers~$p$ with
          \[ 5 \le p \le 1 + \Bigl\lfloor \log_2\Bigl(1 + \frac{\log 4 + \eps(c)/\sqrt{c}}%
                                                               {\log(1 + 1/\sqrt{c})}\Bigr) \Bigr\rfloor ,
          \]
          then all $f^n$ are irreducible.
    \item If $f^2$ is irreducible, $c > 50$, and $\tilde{q}(c) \ge 1.15 c^{-1/30}$,
          then all~$f^n$ are irreducible.
  \end{enumerate}
\end{corollary}

We note that case~(1) gives another proof of case~(4) of Theorem~\ref{main}
for positive~$c$.

For large~$c$, the bound on~$n$ in case~(2) of the corollary is close
to $1 + \log_2(3 + (\sqrt{c} + \tfrac{1}{2}) \log 4)$.

\begin{proof}
  We recall that all~$f^n$ are irreducible when $f^m$ is irreducible for some~$m$
  and $a_n(c)$ or $b_n(c)$ is not a square for all $n > m$.
  \begin{enumerate}
    \item If $c$ is positive and odd, then $f$ is irreducible and $f^2$ is also irreducible
          (since $c$ is not of the form $4m^2(m^2-1)$, compare Proposition~\ref{irredprop}).
          By Proposition~\ref{a3nonsquare}, $a_3(c)$ is never a square when $c > 0$.
          By Proposition~\ref{Prop}, $a_n(c)$ is not a square for all $n \ge 4$, so
          the claim follows.
    \item If $c$ is even and $n > m$, then
          \[ n \ge 1 + \log_2\Bigl(1 + \frac{\eps(n,c)}{\sqrt{c} \log q(c)}
                                     + \frac{\log 4}{\log q(c)}\Bigr) ,
          \]
          since $q(c) > 1 + 1/\sqrt{c}$ and $\eps(n,c) < \eps(c)$. So by Proposition~\ref{Prop},
          $a_n(c)$ is not a square, and the claim again follows.
    \item Let $m$ be as in~(2). Then
          $a_n(c)$ is not a square for $n > m$. For $3 \le n \le m$, $a_n(c)$
          is not a square by assumption (or by  Proposition ~\ref{a3nonsquare} for $n = 3$) if $n$ is prime.
          Otherwise, $n$ is divisible by~$4$ or by an odd prime $p \le m$; then it follows
          that $a_n(c)$ is not a square either, because $(a_n(c))$ is a rigid divisibility
          sequence by Proposition~\ref{divisibility}
          and neither $a_4(c)$ (by Proposition~\ref{a3nonsquare} again) nor $a_p(c)$ is a square.
    \item First note that $2^{2/15} \eps(c)^{1/15} < 1.15$ when $c > 50$.
          The stated inequality then implies that the bound on~$n$ in the second
          statement of Proposition~\ref{Prop} is $< 5$.
    \qedhere
  \end{enumerate}
\end{proof}

We remark that recent work by one of the authors~\cite{stoll2} shows
that $a_5(c)$ is never a square
when $c \neq 0$, which allow us to replace ``$5$'' by~``$7$'' in case~(3)
of the corollary and the condition in case~(4) by ``$\tilde{q}(c) \ge 1.034 c^{-1/126}$''.

We can use case~(4) of Corollary~\ref{C:improve5.4} to deduce case~(5)
of Theorem~\ref{main}; case~(6) of this theorem follows by a similar argument.

\begin{proof}[Proof of cases (5) and (6) of Theorem \ref{main}]
  We can assume that $c > 50$ and $c$ is even, since negative $c$ are dealt
  with by case~(1) and odd $c$ are covered by case~(4)
  of the theorem; the few positive even~$c \le 50$ can be checked individually
  by the methods of this section. Then the
  assumptions of case~(5) imply that $f_r^2$ is irreducible by Proposition~\ref{irredprop}.
  Since when $c$ is a square,
  $c$ cannot be of the form $4m^2(m^2-1)$ either, this is also true in case~(6).

  We first consider case~(5).
  Assume that $c = uv$ with $u$ and~$v$ coprime and (say) $|u|$ a square.
  Then $|v| \ge \prod\limits_{p : 2 \nmid v_p(c)} p^{v_p(c)}$
  and $|u| \le \prod\limits_{p: 2 \mid v_p(c)} p^{v_p(c)}$,
  so that the inequality in the statement implies that $\tilde{q}(c) > 1.15 c^{-1/30}$.
  The claim follows by invoking case~(4) of Corollary~\ref{C:improve5.4}.

  We now consider case~(6). If the claim is false, then there is $n \ge 5$ such
  that $a_n(c)$ is a square. By Lemma~\ref{L:dioph} it follows that we can
  write $c = uv$ with coprime $u$ and~$v$, with $u$ even, such that
  $v^{2^{n-1}-1} - \tfrac{1}{4} u^{2^{n-1}-1} = a_{n-1}(c)$.
  Both $u$ and~$v$ are now squares up to sign, so that we have
  \[ (v^{2^{n-1}-1}, \tfrac{1}{4} u^{2^{n-1}-1}) = \pm (x^2, y^2) \]
  with coprime integers $x$ and~$y$, which implies that
  \begin{equation} \label{E:x2y2}
    x^2 - y^2 = \pm a_{n-1}(c) .
  \end{equation}
  Recall that $a_{n-1}(c) \equiv 1 \bmod c$. Since $c$ is an even square, $x$ is odd,
  and $y$ is even, we obtain the congruence $1 \equiv \pm 1 \bmod 4$, which shows that
  we must have the positive sign in~\eqref{E:x2y2}.
  Let $p \not\equiv 1 \bmod 4$ be a prime dividing~$c$; since $c$ is a square, $p^2 \mid c$.
  It follows that $x^2 \equiv y^2 + 1 \bmod p^2$, and since
  $-1$ is a non-square mod~$p^2$, $p \mid x$ is impossible, so that we must have $p \mid y$.
  This in turn implies that $|u| \ge \prod\limits_{p : p \not\equiv 1 \bmod 4} p^{v_p(c)}$
  and $|v| \le \prod\limits_{p : p \equiv 1 \bmod 4} p^{v_p(c)}$.
  The inequality in the statement then implies that $u/v > 1.15 c^{-1/30}$.
  This contradicts the second inequality in~\eqref{E:uveps}, so that we can
  conclude as in the proof of Proposition~\ref{Prop} that $a_n(c)$ cannot
  be a square, a contradiction.
\end{proof}

\section{A fast algorithm and the proof of case (7) of Theorem \ref{main}} \label{S:algo}

In this section we always assume that $c \ge 4$ is an even integer.
Fix $n \ge 5$ and assume that $a_n = a_n(c)$ is a square.
Set $N = 2^{n-1}-1$.
By Lemma~\ref{L:dioph}, we can write $c = uv$ with $u$ and~$v$ coprime integers
and $u$ even such that
\begin{equation} \label{E:rel}
  v^N - \frac{1}{4} u^N = a_{n-1}(c) .
\end{equation}

We now consider equation~\eqref{E:rel}
as a relation between real numbers. First note that for $c \ge 6$, we have
(using Lemma~\ref{L:anasymp} for the second inequality)
\[ c^{2^{n-2}-1} + c^{2^{n-2}-2}
    \le a_{n-1}(c)
    \le c^{2^{n-2}} F(c)
    \le c^{2^{n-2}-1} + 2 c^{2^{n-2}-2} ,
\]
so~\eqref{E:rel} implies that
\[ v^N - \frac{1}{4} u^N = (uv)^M + \lambda (uv)^{M-1} \]
with $1 \le \lambda \le 2$, where $M = 2^{n-2}-1$ (so that $N = 2M + 1$).

We now set $\theta := 2^{1/N}$, $x := \theta^{-1} u$ and $y := \theta v$; this gives
\[ y^N - x^N = 2 (xy)^M + 2 \lambda (xy)^{M-1} . \]
Writing
\[ z := \frac{(xy)^M}{x^{2M} + x^{2M-1} y + \ldots + y^{2M}} > 0 \]
and recalling that $N = 2M + 1$, this leads to
\begin{equation} \label{E:xy}
  y - x = 2 \Bigl(1 + \frac{\lambda}{xy}\Bigr) z .
\end{equation}
We want to estimate~$z$. We expect that $z$ is close to~$1/N$, which
is the value we obtain when $x = y$. Since $x^{2M-k} y^k + x^k y^{2M-k} \ge 2 (xy)^M$,
it follows that
\[ z \le \frac{1}{N} . \]
Since $xy = uv = c \ge 6$, we see that $y-x$ has to be small:
\begin{equation} \label{E:diffbound}
  0 < y - x < \frac{3}{N} .
\end{equation}

We get a lower bound on~$z$ as follows.
Write $w_k := x^k + x^{k-1} y + \ldots + y^k$.
We consider
\begin{multline}
 \frac{xy}{(y-x)^2} (1 - Nz) 
       =  \frac{xy w_{2M} - N (xy)^{M+1}}{(y-x)^2 w_{2M}} \\
       =  \sum_{j=1}^M \frac{(xy)^{M+1-j} (y^j - x^j)^2}{(y-x)^2 w_{2M}} 
       =  \sum_{j=0}^{M-1} \frac{(xy)^{M-j} w_j^2}{w_{2M}} .
\end{multline}
We note that $(xy)^{M-j} w_j^2$ is given by 
\[ 
     x^{M+j} y^{M-j} + 2 x^{M+j-1} y^{M-j+1} + \ldots + (j+1) x^M y^M + \ldots + x^{M-j} y^{M+j},
\]
which is at most $(j+1) w_{2M}$, and this gives that
\[ \frac{xy}{(y-x)^2} (1 - Nz) \le \sum_{j=0}^{M-1} (j+1) = \frac{M (M+1)}{2}. \]
Thus, using~\eqref{E:diffbound} for the second inequality,
\[ z \ge \frac{1}{N} - \frac{M(M+1)}{2N} \, \frac{(y-x)^2}{xy}
     \ge \frac{1}{N} - \frac{9 M(M+1)}{2 N^3} \, \frac{1}{xy} .
\]
Using this, $1 \le \lambda \le 2$, and $0 < \lambda/(xy) \le 1/3$ in~\eqref{E:xy}, we obtain
\[ \Bigl| y - x - \frac{2}{N} \Bigr| \le \frac{4}{N}\,\frac{1}{xy} . \]
Going back to our original integral variables
$u$ and~$v$, this final bound is equivalent to
\begin{equation} \label{E:bound1}
  \Bigl| \theta^2 v - u - \frac{2 \theta}{N} \Bigr| \le \frac{4 \theta}{N} \, \frac{1}{uv} .
\end{equation}
We want to replace $1/(uv)$ on the right by~$1/v^2$. The following lemma allows
us to do that.

\begin{lemma} \label{L:cuvbound}
  Let $c \ge 4$ be even. We assume that $a_n(c)$ is a square for some $n \ge 2$
  and take $u$ and~$v$ as in~\eqref{E:rel}. Then, with $N = 2^{n-1} - 1$
  and $\theta = 2^{1/N}$,
  \[ (3 - 2 \sqrt{2})^{1/N} < \frac{u}{\theta^2 v} < (3 + 2 \sqrt{2})^{1/N} \]
  and
  \[ (3 - 2 \sqrt{2})^{1/N} < \frac{\theta^2 v}{u} < (3 + 2 \sqrt{2})^{1/N} . \]
  In particular,
  \[ \frac{1}{uv} < \frac{(3 + 2 \sqrt{2})^{1/N}}{\theta^2} \, \frac{1}{v^2}
     \qquad\text{and}\qquad
     c > \theta^2 (3 - 2 \sqrt{2})^{1/N} v^2 .
  \]
\end{lemma}

\begin{proof}
  Note that for $c \ge 4$, we have $a_{n-1}(c) < 2 c^{2^{n-2}-1} \le c^{N/2}$.
  Using this in~\eqref{E:rel} and dividing by~$v^N$, this gives
  \[ \Bigl| 1 - \Bigl(\frac{u}{\theta^2 v}\Bigr)^N \Bigr| < 2 \sqrt{\Bigl(\frac{u}{\theta^2 v}\Bigr)^N} . \]
  Set $\mu := (u/(\theta^2 v))^{N/2} > 0$. Rearranging, we obtain that
  \[ (\mu - 1)^2 < 2 \qquad\text{and}\qquad (\mu + 1)^2 > 2 , \]
  which gives
  \[ (\sqrt{2} - 1)^2 < \mu^2 = \Bigl(\frac{u}{\theta^2 v}\Bigr)^N < (\sqrt{2} + 1)^2 , \]
  from which the bounds in the statement are easily derived.
\end{proof}

\begin{corollary} \label{C:inhombound}
  Let $c \ge 6$ be even. We assume that $a_n(c)$ is a square for some $n \ge 2$
  and take $u$ and~$v$ as in~\eqref{E:rel}. Then, with $N = 2^{n-1} - 1$
  and $\theta = 2^{1/N}$,
  \begin{equation} \label{E:latticebound}
    \Bigl| \theta^2 v - u - \frac{2 \theta}{N} \Bigr|
        < \frac{4 (3 + 2 \sqrt{2})^{1/N}}{\theta N} \, \frac{1}{v^2} .
  \end{equation}
\end{corollary}

\begin{proof}
  This follows immediately from~\eqref{E:bound1} and Lemma~\ref{L:cuvbound}.
\end{proof}

We can use the estimate~\eqref{E:latticebound}
to compute a large lower bound on~$v$ (and therefore on~$c$, by Lemma~\ref{L:cuvbound}),
in the following way. We set
\[ \delta := \frac{4 (3 + 2 \sqrt{2})^{1/N}}{\theta N} . \]
Choose some
$\eps > 0$ (roughly of size $B^{-2}$ when $B$ is the desired lower bound for~$v$).
Let $\Lambda_\eps \subset \R^2$ be the lattice generated by the vectors
$(\eps, \theta^2)$ and~$(0, -1)$. Use lattice basis reduction to find the
minimal squared euclidean distance~$\sigma(\eps)$ between a lattice point and~$(0, 2 \theta/N)$.
Now, assuming that $(v, u) \in \Z^2$
satisfies $|\theta^2 v - u - 2 \theta/N| < \delta/v^2$
(which follows by Corollary~\ref{C:inhombound} for suitable $(v,u)$
if $a_n(c)$ is a square), we see that
\[ \sigma(\eps) \le (\eps v)^2 + |\theta^2 v - u - 2\theta/N|^2 < \eps^2 v^2 + \delta^2 v^{-4} . \]
If the polynomial $\eps^2 X^3 - \sigma(\eps) X^2 + \delta^2$ has two positive
roots $0 < \xi_- \le \xi_+$, then it follows that
\[ |v| > \sqrt{\xi_+} \qquad\text{or}\qquad |v| < \sqrt{\xi_-} . \]
If we already know (from Proposition~\ref{Prop} or a previous application
of the method) that $|v|$ must be larger than~$\sqrt{\xi_-}$, then we
get the new lower bound $|v| > \sqrt{\xi_+}$.

Since the covolume of~$\Lambda_\eps$ is~$\eps$, we expect that
$\sigma(\eps) \approx \eps$.
If $\eps$ is sufficiently smaller than~$\delta^2$, then we get
$\sqrt{\xi_-} \approx \sqrt{\delta}/\sqrt[4]{\eps}$ and $\sqrt{\xi_+} \approx 1/\sqrt{\eps}$.

This gives the following algorithm for checking that $a_n(c)$ can never
be a square when $4 \le c \le \theta^2 (3 - 2 \sqrt{2})^{1/N} B^2$, for a large bound~$B$.

\begin{enumerate}
  \item Use Proposition~\ref{Prop} and Lemma~\ref{L:cuvbound} to determine $B_0$ such that
        $|v| \ge B_0$ in any solution of $a_n(c) = \square$. For example, we can take
        \[ B_0 := \left\lceil \frac{(\sqrt{2}-1)^{1/N}}{\theta}
                     \Bigl(\frac{N}{\log 4} - 3\Bigr) \right\rceil ,
        \]
        where $N = 2^{n-1} - 1$ and $\theta = 2^{1/N}$ as usual.
  \item Repeat the following steps until $B_0 > B$.
        \begin{enumerate}
          \item \label{Step_2a} Set $\eps := \gamma \delta^2 B_0^{-4}$ with some~$\gamma \approx 1$.
          \item Compute $\sigma(\eps)$ and $\xi_-$, $\xi_+$.
          \item If $\xi_- \ge B_0^2$ (or does not exist), increase~$\gamma$
                and go to  Step~\eqref{Step_2a} .
          \item Set $B_0 := \left\lceil \sqrt{\xi_+} \right\rceil$.
        \end{enumerate}
\end{enumerate}

If the algorithm terminates, then this gives a proof that $|v| \ge B$
and therefore (by Lemma~\ref{L:cuvbound}) $c > \theta^2 (3 - 2 \sqrt{2})^{1/N} B^2$ in any
solution of $a_n(c) = \square$.

Since this uses real numbers, it does not yet give a method that can
be implemented on a computer. We need to figure out which precision is necessary.
The lattice basis reduction will essentially compute continued fraction
approximations to~$\theta^2$ with numerators and denominators of size roughly~$B_0^2$.
The resulting reduced lattice basis will have lengths of order~$B_0^{-2}$. The vector
that is closest to $(0, 2\theta/N)$ will then have coefficients of order~$B_0^2$
in terms of this lattice basis. We need the resulting minimal distance to be
computed to an accuracy that is somewhat better than~$B_0^{-2}$. This means
that we need more than $6 \log_2 B_0$ bits of precision. In practice, we
work with an integral lattice obtained by scaling and rounding the basis
given above, as follows. We assume that we have computed $\theta$ to
$> 8 \log_2 B_0$ bits of precision. In the following, $\lfloor \alpha \rceil$
denotes any integer~$a$ such that $|\alpha - a| \le 1$. We can then make
the loop in Step~2 of the algorithm above precise in the following way.

\begin{enumerate}
  \item Set $\gamma := \delta^2$.
  \item \label{Step_b} Let $\Lambda \subset \Z^2$ be the lattice generated by
        $(\lfloor \gamma B_0^4 \rceil, \lfloor \theta^2 B_0^8 \rceil)$ and $(0, -B_0^8)$.
  \item Compute the four points of~$\Lambda$ closest to $(0, \lfloor 2 \theta B_0^8/N \rceil)$;
        call $(v_j, u_j)$ (for $j = 1,2,3,4$) their coefficients with respect to the original
        basis of~$\Lambda$ and set \\
        $(a_j, b_j) := (v_j \lfloor \gamma B_0^4 \rceil,
                        v_j \lfloor \theta^2 B_0^8 \rceil - u_j B_0^8 - \lfloor 2 \theta B_0^8/N \rceil)$.
  \item Set $\sigma := \min\limits_{j=1,2,3,4}
                          \bigl(\max\{0, |a_j| - |v_j|\}^2 + \max \{0, |b_j| - |v_j| - 1\}^2\bigr)$.
  \item Set $h(x) = \lfloor \gamma B_0^4 \rceil^2 x^6 - \sigma x^4 + \lceil \delta^2 B_0^{16} \rceil \in \Z[x]$.
  \item If $h(B_0) \ge 0$, then set $\gamma := 2 \gamma$ and go to  Step~\eqref{Step_b} .
  \item Set $B_0 := \max\{x \in \Z_{> B_0} : h(x) \le 0\} + 1$.
\end{enumerate}

The main point here is that $\sigma/B_0^{16}$ is a lower bound for $\sigma(\eps)$,
where $\eps = \gamma/B_0^4$.

If we denote the successive values taken by~$B_0$ by $B_0, B_1, B_2, \ldots$,
then we expect that $B_{k+1} \approx B_k^2/\delta$. So to reach a given bound~$B$,
we will have to make about $\log\log B$ passes through the loop. The computational
cost of the last pass dominates all others; it is polynomial in~$\log B$.

\begin{example}
  We illustrate how the method works in the case $n = 5$.
  The initial lower bound for~$|v|$ is~$B_0 = 8$. The lattice~$\Lambda$
  has basis $(336, 18401670)$ and~$(0, -16777216)$; the target vector
  is~$(0, 2342757)$. We compute $\sigma = 2373638400$ and find that $h(B_0) < 0$.
  We obtain the new lower bound $B_1 = 145$. In the same way, we find
  the successive lower bounds
  \begin{align*}
    B_2 &= 56956 \\
    B_3 &= 1196488139 \\
    B_4 &= 7319637204404186177 \\
    B_5 &= 41458361126834155279142315082592517830 
  \end{align*}
  and so on.
\end{example}

This allows us to verify Conjecture~\ref{stabconj} for all~$c$ up to
a very large bound~$X$ in reasonable time. We just have to run our algorithm
for all $n = p \ge 5$ prime and the corresponding bound~$B$ for~$|v|$, as long
as the initial bound~$B_0$ (which grows roughly like~$2^p$) is less than~$B$.
Using a straight-forward implementation in MAGMA~\cite{magma}, it took
less than 12~minutes (on the laptop of one of the authors) to prove the following,
which by Proposition~\ref{irredprop} amounts to a proof of case~(7) of Theorem~\ref{main}.

\begin{proposition} \label{P:big}
  Let $c \in \Zplus$ and set $f(x) = x^2 + 1/c$. If $c \le 10^{1000}$
  and the second iterate $f^2$ is irreducible, then all iterates of~$f$
  are irreducible.
\end{proposition}

\begin{proof}
  For $c \in \{1,2,3\}$, this can be checked by considering $a_n(c)$
  modulo $3$, $5$, and~$11$, respectively. So we can assume that $c \ge 4$.
  By Corollary~\ref{C:improve5.4}, it is enough to show that
  $a_p(c)$ is not a square for even~$c$ as in the statement and
  primes
  \[ 5 \le p \le 1 + \Bigl\lfloor \log_2\Bigl(1 + \frac{\log 4 + \eps(c)/\sqrt{c}}%
                                                       {\log(1 + 1/\sqrt{c})}\Bigr) \Bigr\rfloor
             \le 1 + \Bigl\lfloor \log_2\Bigl(3.01 + \frac{\log 4}{\log(1 + 10^{-500})}\Bigr) \Bigr\rfloor,
  \]
  and this last expression is $1662$.
  (For $c \ge 10400$, we use the bound~\eqref{E:c10400}. For smaller~$c$,
  the expression is much smaller than~$1662$.)
  This is a finite computation using the algorithm described above.
\end{proof}

We remark that in the course of executing the algorithm, it was never
necessary to increase the initial value of~$\gamma$.

We have now at last assembled all the ingredients required to prove Theorem \ref{bigcomp}. \label{bigproof}

\begin{proof}[Proof of Theorem \ref{bigcomp}]
  Let $f_r(x) = x^2 + r$ with $r = 1/c$, $c \in \Z \setminus \{0,-1\}$
  and $|c| \leq 10^9$. If $f^2_r(x)$ is irreducible and $c$ is negative or odd,
  then the claim follows from parts (1) and~(4) of Theorem~\ref{main}
  (recall that $f_r(x)$ and so also $f_r^2(x)$ is reducible when $-c$ is a square).
  If $f^2_r(x)$ is irreducible and $c$ is positive and even, Theorem \ref{bigcomp}
  holds by Proposition~\ref{P:big}.

  If $f_r(x)$ or $f_r^2(x)$ is reducible, then the relevant cases of Theorem~\ref{bigcomp}
  follow from Theorem~\ref{maincases}, Corollaries \ref{g2cor} and~\ref{h2cor},
  and Proposition~\ref{h22cor}.
\end{proof}

\section{Applications to the density of primes dividing orbits} \label{density}

In this section, we prove Theorem \ref{densitythm}, which we restate here for the reader's convenience.
\begin{theorem} \label{densitythm2}
Let $c \in \Z$, let $r = 1/c$, suppose that $-c$ and $c+1$ are non-squares in $\Q$, and assume that Conjecture \ref{nosquare} holds for $c$, i.e. that $\frac{a_{n-1} + \sqrt{a_{n}}}{2}$ is not a square in $\Q$ for all $n \geq 3$. Then
for any $t \in \Q$ we have $D(\{\text{$p$ prime : $p$ divides $O_{f_r}(t)$\}}) = 0$.
\end{theorem}
\begin{remark}
Observe that the hypothesis that $\frac{a_{n-1} + \sqrt{a_{n}}}{2}$ not be a square for $n \geq 2$ is strictly weaker than $a_n$ not being a square for $n \geq 2$; in the latter case the conclusion of Theorem \ref{densitythm2} follows immediately from part (2) of \cite[Theorem 1.1]{zdc}. To prove Theorem \ref{densitythm2}, we must apply \cite[Theorem 1.1]{zdc} in a non-trivial way.
\end{remark}
\begin{remark}
When the hypotheses of Theorem \ref{densitythm2} are satisfied, we also obtain certain information on the action of $G_\Q$ on $T_f(0)$ (see p. \pageref{Tdef} for the definition). The index-two subgroup $G_{\Q(\sqrt{-r})}$ acts on both $T_f(\sqrt{-r})$ and $T_f(-\sqrt{r})$. Both of these actions are transitive on each level of the tree, i.e., on $f_r^{-n}(\sqrt{-r})$ (resp. $f_r^{-n}(-\sqrt{-r})$), and the images of the maps $G_{\Q(\sqrt{-r})} \to \text{Sym}(f_r^{-n}(\pm \sqrt{-r})) \cong S_{2^n}$ cannot lie in the alternating subgroup.
\end{remark}

\begin{proof}
Let $K = \Q(\sqrt{-r}),$ so that $f_r = (x + \sqrt{-r})(x - \sqrt{-r})$ over $K$. Let $g_1 = (x + \sqrt{-r})$ and $g_2 = (x - \sqrt{-r})$. To apply part (2) of \cite[Theorem 1.1]{zdc}, we  must show that for $i = 1, 2$, $g_i(f_r^{n-1}(0))$ is a non-square in $K$ for all $n \geq 3$, and also that $-g_i(f_r(0))$ is a non-square in $K$. But
$g_i(f_r^{n-1}(0)) = f_r^{n-1}(0) \pm \sqrt{-r}$.  As in the final part of the proof of Lemma \ref{biglemma2}, $f_r^{n-1}(0) \pm \sqrt{-r}$ is a square in $K$ if and only if
$(f_r^{n-1}(0) \pm \sqrt{f_r^{n}(0)})/2$ is a square in $\Q$, which in turn is equivalent to
$(a_{n-1} + \sqrt{a_{n}})/2$ being a square in $\Q$. But by assumption $(a_{n-1} + \sqrt{a_{n}})/2$ is not a square in $\Q$ for any $n \geq 3$. Moreover, $-g_i(f_r(0)) = -r \mp \sqrt{-r}$, which is a square in $K$ if and only if $(-r \pm \sqrt{r^2+r})/2$ is a square in $\Q$. Because $c+1$ is not a square in $\Q$, it follows that $r^2 + r$ is not a square in $\Q$ either, proving that $-g_i(f_r(0))$ is not a square in $\Q$. 

Therefore we may apply part (2) of \cite[Theorem 1.1]{zdc} twice to show
\begin{equation} \label{denseeq}
0 = \lim_{B \to \infty} \frac{\#\{\p \in S : N({\mathfrak{p}}) \leq B\}}{\#\{{\mathfrak{p}} : N({\mathfrak{p}}) \leq B\}},
\end{equation}
where $N(\p)$ is the norm of the ideal $\p$ and $S$ is the set of primes $\p$ in the ring of integers $\mathcal{O}_K$ of $K$ that divide $g_i(f_r^{n-1}(t))$ for at least one value of $i \in \{1,2\}$ and at least one $n \geq 2$.

If we exclude the finite set of ramified primes, then the primes $\p$ in $\mathcal{O}_K$ come in two flavors: those with norm $p$, where necessarily $p$ splits in $\mathcal{O}_K$; and those with norm $p^2$, where necessarily $p$ is inert in $\mathcal{O}_K$.
Note that $\#\{n \leq B : \text{$n = p^2$ for some prime $p$}\}$ has asymptotic density zero relative to $\#\{n \leq B : \text{$n = p$ for some prime $p$}\}$, and so \eqref{denseeq} is equivalent to
\begin{equation} \label{denseeq2}
0 = \lim_{B \to \infty} \frac{\#\{\p \in S : N({\mathfrak{p}}) = p \leq B\}}{\#\{{\mathfrak{p}} : N({\mathfrak{p}}) = p \leq B\}}.
\end{equation}

Suppose $\p$ in $S$, and say $\p \mid g_i(f_r^{n-1}(t))$ for $n \geq 2$. Then 
$$N(\p) \mid N_{K/\Q}(g_i(f_r^{n-1}(t))) = f_r^{n}(t),$$ 
where $N_{K/\Q}$ is the usual field norm. Let $p = \Z \cap \mathcal{O}_K$ be the prime lying below $\p$. Note that $N(\p) = p$ if $p$ splits in $\mathcal{O}_K$, i.e. if $-r$ is a quadratic residue modulo $p$, and $N(\p) = p^2$ otherwise. But $0 \equiv f_r(f_r^{n-1}(t)) \equiv (f_r^{n-1}(t))^2 + r \bmod{p}$ and hence $-r$ must be a quadratic residue modulo $p$. Thus $N(\p) = p$. It follows that the numerator of \eqref{denseeq2} is $2\#\{p : \text{$p \leq B$ and $p$ divides $O_f(t)$}\}$. Clearly the denominator is
$$2\#\{p : \text{$p \leq B$ and $-r$ is a quadratic residue modulo $p$}\}.$$
 But by quadratic reciprocity and Dirichlet's theorem on primes in arithmetic progressions, the latter is asymptotic to $\#\{p : p \leq B\}$. It follows that $D(\{p : \text{$p$ divides $O_f(t)$}\}) = 0$, as desired.
\end{proof}

\bibliographystyle{plain}

\end{document}